\documentclass[12pt]{article}
\usepackage{latexsym}
\usepackage{graphicx}
\usepackage{subcaption}
\usepackage{amsmath}
\usepackage{amsfonts}
\usepackage{dsfont}
\usepackage{amssymb}
\usepackage{amsthm}
\usepackage{enumerate}
\usepackage[xcdraw,table]{xcolor}

\usepackage{algorithmicx}
\usepackage{algcompatible}
\usepackage{algpseudocode}
\usepackage[table]{xcolor}

\usepackage{appendix}
\newtheorem{theorem}{Theorem}[section]
\newtheorem{remark}{Remark}[section]
\newtheorem{lemma}{Lemma}[theorem]

\newtheorem{assumption}{Assumption}[section]

\newcounter{probnum}
\definecolor{tabblue}{rgb}{.870588,.905882,.94902}
\definecolor{gray}{rgb}{0.7,0.7,0.7}
\definecolor{black}{rgb}{0,0,0}
\definecolor{white}{rgb}{1,1,1}
\definecolor{blue}{rgb}{0.0,0.0,1}

\definecolor{green}{rgb}{0,0.5,0}

\definecolor{yellow}{rgb}{1,0.549,0}

\definecolor{red}{rgb}{0.6,0.0,0.0}

\definecolor{darkred}{rgb}{0.9,0.4,0}

\definecolor{purple}{rgb}{0.58,0,0.827}

\definecolor{backgcode}{rgb}{0.97,0.97,0.8}
\definecolor{Brown}{cmyk}{0,0.81,1,0.60}
\definecolor{OliveGreen}{cmyk}{0.64,0,0.95,0.40}
\definecolor{CadetBlue}{cmyk}{0.62,0.57,0.23,0}

\DeclareMathOperator*{\argmin}{arg\,min~}

\newcommand{\qu}[1]{``{#1}''}
\newcommand{\bv}[1]{\boldsymbol{#1}}

\newcommand{\rcontrol}{\tilde{r}}

\newcommand{\bSigma}{\bv{\Sigma}}
\newcommand{\bSigmaw}{\bv{\Sigma}_{W}}
\newcommand{\bSigmaZ}{\bv{\Sigma}_{Z}}

\newcommand{\sigmawij}[2]{\sigma_{#1\,,#2}}
\newcommand{\sigmazi}[1]{\rho_{#1}}

\newcommand{\gammaZ}{\bv{\gamma}_Z}

\newcommand{\curlyV}{\mathcal{V}}

\newcommand{\muvec}{\bv{\mu}}

\newcommand{\ybar}{\bar{y}}
\newcommand{\Ybar}{\bar{Y}}
\newcommand{\YbarT}{\Ybar_T}
\newcommand{\YbarC}{\Ybar_C}
\newcommand{\ybarT}{\ybar_T}
\newcommand{\ybarC}{\ybar_C}

\newcommand{\half}{\frac{1}{2}}

\newcommand{\B}{\bv{B}}

\renewcommand{\r}{\bv{r}}

\newcommand{\Z}{\bv{Z}}

\newcommand{\I}{\bv{I}}
\newcommand{\J}{\bv{J}}
\newcommand{\Y}{\bv{Y}}
\newcommand{\W}{\bv{W}}

\newcommand{\x}{\bv{x}}

\newcommand{\w}{\bv{w}}

\newcommand{\tauhat}{\hat{\tau}}

\newcommand{\onevec}{\bv{1}}
\newcommand{\zerovec}{\bv{0}}

\newcommand{\y}{\bv{y}}

\renewcommand{\v}{\bv{v}}
\renewcommand{\muvec}{\bv{\mu}}

\newcommand{\z}{\bv{z}}

\newcommand{\bbeta}{\bv{\beta}}

\newcommand{\reals}{\mathbb{R}}

\newcommand{\naturals}{\mathbb{N}}

\newcommand{\beqn}{\vspace{-0.25cm}\begin{eqnarray*}}
\newcommand{\eeqn}{\end{eqnarray*}}
\newcommand{\bneqn}{\vspace{-0.25cm}\begin{eqnarray}}
\newcommand{\eneqn}{\end{eqnarray}}
\newcommand{\benum}{\begin{enumerate}}
\newcommand{\eenum}{\end{enumerate}}

\newcommand{\parens}[1]{\left(#1\right)}

\newcommand{\prob}[1]{\mathbb{P}\parens{#1}}

\newcommand{\bracks}[1]{\left[#1\right]}
\newcommand{\braces}[1]{\left\{#1\right\}}

\newcommand{\abss}[1]{\left|#1\right|}
\newcommand{\norm}[1]{\left|\left|#1\right|\right|}
\newcommand{\normsq}[1]{\norm{#1}^2}

\newcommand{\inverse}[1]{\parens{#1}^{-1}}

\newcommand{\diag}[1]{\text{Diag}\bracks{#1}}

\newcommand{\tr}[1]{\text{tr}\bracks{#1}}

\newcommand{\expe}[1]{\mathbb{E}\bracks{#1}}

\newcommand{\expesub}[2]{\mathbb{E}_{\,#1}\bracks{#2}}

\newcommand{\var}[1]{\mathbb{V}\text{ar}\bracks{#1}}

\newcommand{\varsub}[2]{\mathbb{V}\text{ar}_{#1}\bracks{#2}}

\newcommand{\msesub}[2]{\mathbb{M}\text{SE}_{#1}\bracks{#2}}

\newcommand{\cov}[2]{\mathbb{C}\text{ov}\bracks{#1,\,#2}}

\newcommand{\oneover}[1]{\frac{1}{#1}}
\newcommand{\overtwo}[1]{\frac{#1}{2}}

\newcommand{\bernoulli}[1]{\mathrm{Bernoulli}\parens{#1}}
\newcommand{\betanot}[2]{\mathrm{Beta}\parens{#1,\,#2}}

\newcommand{\normnot}[2]{\mathcal{N}\parens{#1,\,#2}}

\newcommand{\poisson}[1]{\mathrm{Poisson}\parens{#1}}

\newcommand{\weibullnot}[2]{\mathrm{Weibull}\parens{#1,\,#2}}

\newcommand{\brho}{\boldsymbol{\rho}}

\allowdisplaybreaks 
\usepackage{authblk}
\usepackage{xcolor}

\usepackage{natbib}
\usepackage{multirow}

\usepackage[margin=1.02in]{geometry}

\newcommand{\ourtitle}{The Optimality of Blocking Designs in Equally and Unequally Allocated Randomized Experiments with General Response}
\title{
\ourtitle
}

\author[1]{David Azriel\thanks{Electronic address: \texttt{davidazr@technion.ac.il}}}
\author[2]{Abba M. Krieger\thanks{Electronic address: \texttt{krieger@wharton.upenn.edu}}}
\author[3]{Adam Kapelner\thanks{Electronic address: \texttt{kapelner@qc.cuny.edu}; Principal Corresponding author}}

\affil[1]{\small Faculty of Industrial Engineering and Management, The Technion, Haifa, Israel}
\affil[2]{\small Department of Statistics, The Wharton School of the University of Pennsylvania, USA}
\affil[3]{\small Department of Mathematics, Queens College, CUNY, USA}

\begin{document}
\maketitle
  
\begin{abstract}
We consider the performance of the difference-in-means estimator in a two-arm randomized experiment under common experimental endpoints such as continuous (regression), incidence, proportion and survival. We examine performance under both equal and unequal allocation to treatment groups and we consider both the Neyman randomization model and the population model. We show that in the Neyman model, where the only source of randomness is the treatment manipulation, there is no free lunch: complete randomization is minimax for the estimator's mean squared error. In the population model, where each subject experiences response noise with zero mean, the optimal design is the deterministic perfect-balance allocation. However, this allocation is generally NP-hard to compute and moreover, depends on unknown response parameters. When considering the tail criterion of Kapelner et al. (2021), we show the optimal design is less random than complete randomization and more random than the deterministic perfect-balance allocation. We prove that Fisher's blocking design provides the asymptotically optimal degree of experimental randomness. Theoretical results are supported by simulations in all considered experimental settings.\\~\\
{\bf Keywords:} experimental design, optimal design, blocking, unequal randomization, restricted randomization
\vfill
\end{abstract}
 
\maketitle

\section{Background}\label{sec:intro}

We consider a classic problem: an experiment with $2n$ \emph{subjects} (\emph{individuals}, \emph{participants} or \emph{units}). The experiment has two \emph{arms} (\emph{treatments}, \emph{manipulations} or \emph{groups}), which we will call treatment and control denoted $T$ and $C$, and the number of treatments denoted $n_T$ and number of controls denoted $n_C$. We consider one clearly defined \emph{outcome} (\emph{response} or \emph{endpoint}) of interest for the $2n$ subjects denoted $\y = \bracks{y_1, \ldots, y_{2n}}^\top$, which is of general type: either continuous, incidence, count or uncensored survival. Subjects have $p$ observed subject-specific \emph{covariates} (\emph{measurements}, \emph{features} or \emph{characteristics}) denoted $\x_i \in \reals^p$ for the $i$th subject. The setting we investigate is where all $\x_i$'s are known beforehand and considered fixed. This non-sequential setting was studied by \citet{Fisher1925} when assigning treatments to agricultural plots and is still of great importance today. In fact, this setting is found in \qu{many phase I studies [that] use `banks' of healthy volunteers ... [and] ... in most cluster randomised trials, the clusters are identified before treatment is started} \citep[page 1440]{Senn2013}. The more common \emph{sequential} experimental setting where subjects arrive one-by-one and must be assigned on the spot (hence all $\x_i$'s are not seen beforehand), we will leave for future work. \citet{Sverdlov2019} provides an extensive discussion and simulation of common sequential unequal allocation designs.

Formally defined, the \emph{randomization} (\emph{allocation} or an \emph{assignment}) is a vector $\w = \bracks{w_1, \ldots, w_{2n}}^\top$ whose entries indicate whether the subject received $T$ (coded numerically as +1) or $C$ (coded numerically as -1) and thus this vector is an element in $\braces{-1,+1}^{2n}$. After the sample is provided, the only degree of control the experimenter has is to choose the entries of $\w$. The process that results in the choices of such allocations of the $2n$ subjects to the two arms is termed an experimental \emph{design} (\emph{strategy}, \emph{algorithm}, \emph{method} or \emph{procedure}). Experimental design is thus a generalized multivariate Bernoulli random variable that we denote $\W$ and its number of allocations as $N_{\W}$ where the design can technically be degenerate (i.e., $N_{\W} = 1$). 

We show herein that the optimal asymptotic design for general response is one that blends pure randomization with deterministic optimization. Further, this optimal design is satisfied by Fisher's historic blocking design \textit{as commonly implemented today}. Before we provide more details about our contribution, we define the equal allocation setting and unequal allocation setting (as our result applies to both settings) in the next two subsections. To contextualize our result, we also provide some advantages and disadvantages of both classic and modern experimental designs.

\subsection{Equal Allocation Experimentation}\label{subsec:equal_allocation_intro}

If $n_T = n_C$, the design is termed \textit{equal allocation}. The default design here is the balanced complete randomization design (BCRD, \citep[p. 1171]{Wu1981}). This design has one constraint: $\sum_i w_i = 0$ which is satisfied by $\binom{2n}{n}$ vectors. These vectors are then sampled uniformly. BCRD has been given the high distinction of \qu{the gold standard}. \citet[p. 245]{Cornfield1959} gives two reasons: (1) known and unknown covariate differences among the groups will be small and bounded and (2) it is the \qu{reasoned basis} for causal inference, a paraphrasing of \citet[p. 14]{Fisher1935}. Any other design (besides the \emph{Bernoulli Trial} design discussed in Section~\ref{subsec:worst_case_criterion}) we will call \qu{restricted} as it is restricted to use fewer assignment vectors. 

Why employ a restricted design? There is a large problem with BCRD that was identified immediately by \citet{Fisher1925}: under some unlucky $\w$'s there are large differences in the distribution of observed covariates between the two groups. The amount of covariate value heterogeneity between groups we term \emph{observed imbalance}. This observed imbalance leads to higher variance during estimation (as defined at the end of Section~\ref{sec:methods}). Through an abuse of terminology, many times the literature denotes this as simply \emph{imbalance}, but this is an ambiguous term since it usually ignores the state of imbalance in the unobserved covariates which also leads to higher variance during estimation  (again, see the end of Section~\ref{sec:methods}).

Mitigating the chance of large observed imbalances is the predominant reason for employing a restricted randomization procedure. These procedures take the observed covariate values as input and compute assignments with small observed imbalance as output. \citet[p. 251]{Fisher1925} wrote \qu{it is still possible to eliminate much of the \ldots heterogeneity, and so increase the accuracy of our [estimator], by laying restrictions on the order in which the strips are arranged}. Here, he introduced the \emph{block design}, a restricted design still popular today, especially in clinical trials. Another early mitigation approach can be found in \citet[p. 366]{Student1938} who wrote that after an unlucky, highly imbalanced randomization, \qu{it would be pedantic to [run the experiment] with [an assignment] known beforehand to be likely to lead to a misleading conclusion}. His solution is for \qu{common sense [to prevail] and chance [be] invoked a second time}. In foregoing the first assignment and \emph{rerandomizing} to find a better assignment, all allocations worse than a predetermined threshold of observed imbalance are eliminated. This classic strategy has been rigorously investigated only recently \citep{Morgan2012, Li2018} in the case of equal allocation. Another idea is to allocate treatment and control among similar subjects by using a pairwise matching (PM) design \citep{Greevy2004}. This design creates $n$ pairs (i.e., blocks with size 2) and randomizes the two subjects to T/C or C/T. All of the above-mentioned restricted designs (blocking, rerandomization, PM) are still \qu{random} in the sense that there are still exponentially many possible allocations.

One may wonder about the \qu{optimal} design strategy. Regardless of how optimality is defined and regardless of whichever criterion is used to measure performance, we believe this is an impossible question to answer. The multivariate Bernoulli random variable has $2^{2n} - 1$ parameters, an exponentially large number \citep[Section 2.3]{Teugels1990}. Finding the optimal design is thus tantamount to solving for exponentially many parameters using only the $2n$ observations, a hopelessly unidentifiable task. Thus, to find an \qu{optimal} design, the space of designs must be limited. 

One way of limiting the space of designs is to do away with any notion of randomization whatsoever and instead locate a \emph{deterministic optimal design} with one optimal allocation vector $\w_*$, a design which we will call \emph{perfect balance} (PB). An early advocate of this approach was \citet{Smith1918} and a good review of classic works advocating this approach such as \citet{Kiefer1959} is given in \citet[Chapter 3]{Steinberg1984}. A particularly edgy advocate is \citet{Harville1975} who penned an article titled \qu{Experimental Randomization: Who Needs It?}. Although unpalatable to modern ears that have been conditioned to the concept that \qu{randomization is the gold standard}, this work and our previous work showed that under certain assumptions, there is no need for randomness in the design when the inferential focus is on effect estimation (hypothesis testing is a different discussion).

How does one find this unique allocation $\w_*$? One needs to first define a numerical metric of imbalance as a function of the $\x$'s and $\w$. Then, the optimization problem becomes a binary integer programming problem where the unequal allocation requirement would merely alter one of the inputted constraints. These types of problems are usually NP-hard but can be approximated via numerical optimization \citep{Bertsimas2015, Kallus2018}. Sometimes there are polynomial-time algorithms that provide \qu{good-enough} solutions for practical purposes, e.g. branch and bound \citep{Land1960}. Selection of the imbalance metric is arbitrary but recent work by \citet{Kallus2018} posits an appropriate metric to use based on assumptions about the structure of the response function.

There are even other designs that attempt to retain the best of both worlds i.e. randomizing while simultaneously optimizing. For instance, \citet{Krieger2019} starts with BCRD and greedily switches pairs of treatment and control subjects until reaching a locally optimal minimum of observed covariate imbalance. This procedure provably retains most of the randomness of BCRD while achieving provably minuscule observed imbalance. \citet{Zhu2021} also use pair-switching but do so non-greedily and thus can explore the space of forced-balanced allocations more completely. \citet{Harshaw2019} use the Gram-Schmidt Walk design to provide allocations that have high-degree of randomness and low observed imbalance. Another approach is the Finite Selection Model procedure of \citet{Chattopadhyay2022}. They employ a sequentially controlled Markovian random sampling to assign $w_i$'s based on a $D$-optimality criterion. 

\subsection{Unequal Allocation Experimentation}\label{subsec:unequal_allocation_intro}

There are many practical situations where subjects are allocated to the two arms unequally. For example, to reduce the cost of the trial (e.g the treatment is an expensive surgery and the control is an inexpensive placebo), avoid loss of power from drop-out or cross-over, various ethical concerns, for gaining additional information on one manipulation preferentially \citep{Dumville2006}, as well as if the response variance is known to be heterogeneous across manipulations \citep{Wong2008}. \citet{Torgerson1997} makes the general recommendation that unequal allocation should be employed more often in the case of economic concerns. The literature on unequal allocation is more limited than on equal allocation. 

We will now discuss the historic designs above in the context of unequal allocation beginning with BCRD. It is clear BCRD must be altered as it is by definition an exclusively equal allocation design. The analogue design of BCRD will be the unique design where all $\binom{2n}{n_T} = \binom{2n}{n_C}$ vectors satisfy the imbalance constraint $\sum_{i=1}^n w_i = n_T - n_C < 0$. We will term this modified design as the \textit{imbalanced complete randomization design} (iBCRD) although it has other names \citep[Section 3.1.1]{Sverdlov2019}. Since BCRD is a special case of iBCRD, going forward we will refer to both as i/BCRD. Blocking in the case of unequal randomization requires each blocks' allocations to largely respect the allocation ratio; and this design is usually termed the \emph{permuted block design} \citep[Section 2.2.2]{Ryeznik2018}. As for PM, this design must be modified as it is by definition an equal allocation design. We talk about this modification in Section~\ref{sec:discussion}. 

Most of the extant unequal allocation literature focuses on (a) computing the allocation proportions and (b) how to then allocate sequentially. The first attempt at (a) is based on heterogeneous treatment costs: \citet{Nam1973} assumes a continuous endpoint distributed normally. \citet{Tymofyeyev2007} propose an optimization framework for (a) in the binary response case and shows it is usually convex. \citet{Feng2018} employ this framework for (b). \citet{Sverdlov2019} propose a procedure for (a) for two or more treatment groups with a normally distributed regression endpoint. Their objective is to minimize the overall cost of the experiment subject to a minimum requirement for statistical power. Further, they discuss performance under both the \emph{population model} and \emph{randomization model}, two different worldviews that are fundamental to the differing results in this work. We view this work as complementary; the extensive work on (a) can be applied directly herein as a pre-step to choose the allocation ratio before implementing the design we recommend. 

We are not aware of work that discusses unequal-allocation design vis-a-vis the consideration of minimizing observed imbalance in the non-sequential setting where all $\x$'s are known a priori. 

\subsection{The Randomness-Imbalance Tradeoff}

The above gives a short introduction to the wide range of design ideas. At the \qu{extremes} of this range, there is i/BCRD on one side and the minimized observed imbalance deterministic PB design on the other side. Everything in between these extremes attempts to compromise between observed imbalance and degree of randomness. 

Is randomness worth sacrificing for less observed imbalance? This question was investigated in \citet{Kapelner2021} for the case of a continuous outcome under equal treatment-control allocation (henceforth referred to as \qu{our previous work}). There was a threefold answer: (I) if the experimenter is concerned that nature may provide the worst possible unobserved covariates then one should never sacrifice any randomness, (II) if the experimenter treats the unobserved covariates as random noise that can be averaged over, then one should sacrifice all randomness and (III) if the experimenter is concerned of tail events of unobserved covariates then one should sacrifice some randomness. A related work is \citet{Nordin2022}, which studies a measure of how correlated different assignment vectors in a design are with each other. They show that the variance of the MSE increases linearly in this measure. Therefore, in order to avoid high variance of the MSE, one needs to consider designs where this correlation metric is small.

The purpose of this work is to reconsider the above question for the most common endpoints in experiments and clinical trials: continuous, incidence, survival, count and proportion and for treatment allocation that are both equal and unequal across treatments. We prove that our previous work's results generalize. We further uncover the asymptotic optimality of blocking designs in (III). We define our experimental setting formally in Section~\ref{sec:setup}, discuss criterions of performance and theoretical results in Section~\ref{sec:methods}, provide simulation evidence of our results in Section~\ref{sec:simulations} and conclude in Section~\ref{sec:discussion}.


\section{Problem Setup}\label{sec:setup}

Consider the potential outcome framework of \citet{Rubin2005} and let $\y_T$ and $\y_C$ denote the vectors of all subject responses under the treatment and control respectively. For each $i=1, \ldots, 2n$ the observed response $y_i$ is $y_{T,i}$ (respectively, $y_{C,i}$) if $w_i=+1$ (respectively, $w_i=-1$). Let $\W$ denote the possible unequal-allocation design which produces vectors where $n_T \leq n_C$ with $n_T$ and $n_C$ fixed. In a vector form, the $2n$ random responses $\Y$ can be expressed as

\bneqn\label{eq:potential_outcomes}
\Y = \half\Big(\y_T + \y_C + \diag{\W} (\y_T - \y_C)\Big).
\eneqn

\noindent under the Stable Unit Treatment Values Assumption. Note that the randomness in the responses in Equation~\ref{eq:potential_outcomes} is only due to the randomness in the design $\W$. This limitation of randomness to the design exclusively is known as the \emph{randomization model} or \emph{Neyman Model} \citep[Chapter 6.3]{Freedman2008, Lin2013, Rosenberger2016} in contrast to a \emph{population model}, which we will describe in Section~\ref{subsec:mean_criterion}. 

We consider the sample average treatment effect as our estimand,

\bneqn\label{eq:estimand}
\tau := \oneover{2n} \onevec_{2n}^\top(\y_T - \y_C).
\eneqn

\noindent This estimand is a natural target of inference for all outcome types. In the incidence setting, this estimand is the mean probability difference (also called the \emph{risk difference}). Sometimes, practitioners prefer the risk ratio, odds ratio or log odds ratio. Future work can examine additional estimands beyond the mean probability difference in the setting of the binary outcome.

An intuitive estimate for our estimand is $\tauhat = \ybarT - \ybarC$ where $\ybarT$ and $\ybarC$ are the average response in the treatment and control groups respectively when computed for a specific $\w$ from any design $\W$. Using Equation~\ref{eq:potential_outcomes} and the elementary computation found in Section~\ref{app:estimator_computation} of the Supplementary Material, we can conveniently express this estimate's estimator as 

\bneqn\label{eq:estimator}
\tauhat = \oneover{2n} \parens{ 
    \onevec_{2n}^\top \parens{\frac{
            \y_T
        }{r} - 
        \frac{
            \y_C
        }{\rcontrol}
    }+
    \W^\top \parens{
        \frac{
            \y_T
        }{r} + 
        \frac{
            \y_C
        }{\rcontrol}
    }
}.
\eneqn

\noindent where $r := n_T / n \in (0, 1]$ and $\rcontrol := n_C / n = 2 - r \in [1, 2)$. This notation has the benefit that substitution of $r = \rcontrol = 1$ (or equivalently, $r - \rcontrol = 0$) implies the equal allocation setting (the more common experimental setting) providing for easy visual simplification of our expressions going forward.

In order to prove results about the performance of this estimator under different designs, we now make two assumption about the design space: 

\begin{assumption}\label{ass:constrained_balance}
All assignments respect a specific equal or unequal allocation, i.e. $\w^\top \onevec_{2n} = (r - \rcontrol)n$ for all $\w$ in the support of $\W$. Note this quantity is zero under equal allocation.
\end{assumption}

\begin{assumption}\label{ass:equal_chance_T}
All subjects have equal probability of being assigned the treatment condition, i.e. $\expe{\W} = \half(r - \rcontrol)\onevec_{2n}$. Note this quantity is the zero vector under equal allocation.
\end{assumption}

\noindent Computation found in Section~\ref{app:estimator_unbiasedness} of the Supplementary Material shows our estimator is unbiased under Assumption~\ref{ass:equal_chance_T}. The mean squared error (MSE) over all assignments is thus its variance which can be expressed as the quadratic form

\bneqn\label{eq:mse}
\msesub{\W}{\tauhat} = \oneover{4n^2} \v^\top \bSigmaw \v
\eneqn

\noindent where the defining matrix $\bSigmaw$ denotes $\var{\W}$, the variance-covariance matrix of the design, and $\v := \y_T/r + \y_C/\rcontrol$, the \qu{potential outcome vector} \citep[page 6]{Harshaw2019}.

We will be analyzing performance for many types of responses. It will be useful to express potential outcome as a sum of two components: (a) $\mu_{T,i}, \mu_{C,i}$ which are a function of the observed covariates $\x_i$ (b) $z_{T,i}, z_{C,i}$ that which is a function of the unobserved covariates. 
Letting $\muvec_T$ and $\muvec_C$ (and $\z_T$ and $\z_C$) denote the stacked vectors of the observed (and unobserved) covariates' contribution to the responses when $\w = +\onevec_{2n}$ and $\w = -\onevec_{2n}$ respectively, we have

\bneqn\label{eq:simple_components}
\y_T = \muvec_T + \z_T \in \mathcal{Y}, ~~ \y_C = \muvec_C + \z_C \in \mathcal{Y}.
\eneqn

\noindent The response types we consider are formally listed in Table~\ref{tab:simple_models} (column 1). The potential response values belong to different sets $\mathcal{Y}$ depending on response type (column 2). This table also lists common generalized linear models (GLMs) for each of the responses' means (column 4). These means likewise belong to different sets $\Theta$, where $\Theta$ is the space of the $\mu_{T,i}, \mu_{C,i}$'s (column 5). In the population model, the $z_{T,i}, z_{C,i}$'s are conceptualized as randomly-realized \emph{noise} or \emph{errors}. Regardless of what one believes about their origination, the critical fact is they cannot be modeled with the information the experimenter has at hand. 

\begin{table}[ht]
    \centering
    \small
    \begin{tabular}{cc|cc|c}
       Response type name & $\mathcal{Y}$ & Model Name & Common $\mu_{T,i}$ and $\mu_{C,i}$ & $\Theta$ \\ \hline
       
       Continuous & $\reals$ & Linear & $\beta_0 + \bbeta^\top\x_i + \beta_T w_i$ & $\reals$\\ \hline
       
       Incidence & $\braces{0,1}$ &  Log Odds & \multirow{2}{*}{$\inverse{1 + e^{-(\beta_0 + \bbeta^\top\x_i + \beta_T w_i)}}$} & \multirow{2}{*}{$(0, 1)$} \\
       
       Proportion & $(0, 1)$ & Linear & & \\ \hline
       
       Count & $\braces{0,1, \ldots}$ & \multirow{2}{*}{Log Linear} & \multirow{2}{*}{$e^{(\beta_0 + \bbeta^\top\x_i + \beta_T w_i)}$} & \multirow{2}{*}{$(0, \infty)$}\\
       
       Survival & $(0, \infty)$ &  &\\ \hline
       
    \end{tabular}
       
       
       
       
       
       
    \caption{Response types considered in this work. For the example mean response model, the $\beta$'s are unknown parameters to be estimated of which $\beta_T$ is most often of paramount interest to the experimenter.}
    \label{tab:simple_models}
\end{table}

\section{Analyzing Designs}\label{sec:methods}

It is impossible to employ the MSE as the criterion in which we compare different designs since $\z_T$ and $\z_C$ are not observed. Thus, as in our previous work, we do not employ MSE as the performance criterion \textit{directly} but instead remove the dependence on $\z_T, \z_C$ in three different ways to arrive at three different practical criterions (a) the worst-case MSE, (b) the mean MSE and (c) the tail MSE. These are explored in the next three sections respectively.

Before we do so, it is worthwhile to have intuition about what drives the MSE of our estimator. Using the partitioning of the response in Equation~\ref{eq:simple_components} and letting let  $\muvec := \muvec_T  / r + \muvec_C  / \rcontrol$ and $\z := \z_T / r + \z_C / \rcontrol$, we can rewrite Equation~\ref{eq:mse} dropping the constant multiple as 

\beqn
\msesub{\W}{\tauhat} \propto \muvec^\top \bSigmaw \muvec + 2\muvec^\top \bSigmaw \z + \z^\top \bSigmaw \z.
\eeqn

\noindent We will see below that the first term represents imbalance in the observed covariates and the last term represents imbalance in the unobserved covariates.

The matrix $\bSigmaw$ is defined as $\expe{\W\W^\top} - \expe{\W}\expe{\W}^\top$ where the latter term is constant for all designs by Assumption~\ref{ass:equal_chance_T}. Its first term can be expressed as

\bneqn\label{eq:expe_w_w_transpose}
\expe{\W\W^\top} = \sum_{k=1}^{N_{\W}} p_k \w_k \w_k^\top
\eneqn

\noindent where $p_k := \prob{\W = \w_k}$, i.e., the probability that the $k$th allocation is used in the experiment. Thus, the MSE can be alternatively expressed as

\bneqn\label{eq:mse_decomposition_three_terms}
\msesub{\W}{\tauhat} \propto c_{r,n} + \sum_{k=1}^{N_{\W}}p_k \parens{
(\muvec^\top \w_k)^2 + 2\muvec^\top\w_k \z^\top \w_k + (\z^\top \w_k)^2
}
\eneqn

\noindent where $c_{r,n}$ is independent of the design. The first term within the parentheses represents the mean (over all allocations) of imbalance among the observed covariates' mean contribution and the last term within the parentheses represents the mean (over all allocations) of imbalance among the unobserved covariates. In the case of continuous response, linear mean model and $p=1$, the first term becomes $(\x^\top \w_k)^2$ which is the mean difference of the observed covariate values between treatment and control arms. The outer sum with $p_k$ takes the mean of this mean over all allocations properly weighted.

\subsection{The Worst Case MSE}\label{subsec:worst_case_criterion}

We attempt to find the best design under an adversarial nature that provides us the worst possible $\z_T, \z_C$ values,

\beqn
\W_* = \argmin_{\W \in \mathcal{W}} \sup_{\z_T, \z_C} \msesub{\W}{\tauhat}
\eeqn

\noindent where $\mathcal{W}$ denotes the space of all legal designs whose allocations satisfy Assumptions~\ref{ass:constrained_balance} and \ref{ass:equal_chance_T}. In the supremum operation, the space of legal $\z_T$'s and $\z_C$'s will then depend on the response type. Since the $\z$ vectors sum with $\muvec$ to yield the response (Equation~\ref{eq:simple_components}) and the response averages are subtracted to compute $\v$, this inner problem is equivalent to $\sup_{\v \in \curlyV} \v^\top \bSigmaw \v$, where $\curlyV$ is the space of all possible $\v$'s. When $\curlyV$ is unbounded the problem is ill-defined as the inner supremum is infinite. Thus, to make the problem realistic (response metrics in real-world experiments are bounded), we limit the space of vectors to $\curlyV_M \subset \curlyV$. The definition of the subspace $\curlyV_M$ depends on the response. For the continuous case it is natural to bound the space by the Euclidean norm of $\v$ as in \citet{Efron1971} who employs this norm to define ``accidental bias'' bounded at one (without loss of generality). For the count and survival response it is appropriate to bound the sup norm of $\v$. Thus, we define

\bneqn \label{eq:D_M}
\curlyV_M := \left\{ \begin{array}{cc}  \{ \v \in \mathbb{R}^{2n}  : \| \v \| \le 1\} & \text{continuous response}\\
\braces{0, 1/r, 1/ \tilde{r},1/r + 1/ \tilde{r}}^{2n} & \text{incidence}\\
\parens{0, 1/r + 1/ \tilde{r}}^{2n} &  \text{proportion responses}\\
\braces{ \v \in \braces{i/r + j/ \tilde{r} \,:\, i,j \in \naturals_0 }^{n} : \max(\v) \le M } & \text{count}\\
\braces{ \v \in (0,\infty)^{n} : \max(\v) \le M }  &  \text{survival responses}
\end{array}\right.
\eneqn

\noindent where $M$ is a large positive constant such that $\prob{Y_i > M}$ is negligible for all subjects. The design location problem then becomes

\bneqn\label{eq:worst_case_problem}
\W_* = \argmin_{\W \in \mathcal{W}} \max_{\v \in \curlyV_M} \v^\top \bSigmaw \v,
\eneqn

\noindent which is exactly tractable for the case of the response being continuous. In this case, the inner maximum is the largest eigenvalue of $\bSigmaw$. The $\W$ that minimizes this largest eigenvalue is such where each subject is assigned on the basis of a fair coin flip and is either called the \emph{complete randomization} design \citep[Chapter 3.2]{Rosenberger2016} or the \emph{the Bernoulli Trial} \citep[Chapter 4.2]{Imbens2015}. In this design, $\bSigmaw = \I_{2n}$ which mirrors the minimax result found in \citet[Section 2.1]{Kallus2018}. However, this design violates our space of considered designs which are only those that satisfy a prespecified allocation to the two treatments (Assumption~\ref{ass:constrained_balance}). Within the space of allowable designs, we have the following no free lunch result which follows from \citet[Theorem 2.1]{Kapelner2021} and proven in Section~\ref{app:ibcrd_minimax} in the Supplemental Material.

\begin{theorem}\label{thm:iBCRD_minimax_for_continuous}
Under Assumptions \ref{ass:constrained_balance} and \ref{ass:equal_chance_T} and continuous response, i/BCRD is the minimax optimal design in the case of unequal / equal allocation for the worst case MSE criterion (Equation~\ref{eq:worst_case_problem}).
\end{theorem}

\noindent We can think of this result as minimizing the effect of the third term in the MSE decomposition of Equation~\ref{eq:mse_decomposition_three_terms} while not addressing the first term. 

If the response is non-continuous, the max problem embedded within Equation~\ref{eq:worst_case_problem} is generally intractable. In order to make this problem tractable, we further limit the set of designs to block designs.

\subsubsection{The Optimal Design Among Block Designs}\label{subsubsec:block_designs_worst_case}

We denote block designs with $B$ blocks as $\W \in BL(B)$. Each block is of size $n_B$, where $n_B:= 2n / B$, and in each block there will be exactly $\frac{n_T}{2n} n_B$ subjects assigned to the treatment manipulation and $ \frac{n_C}{2n} n_B$ subjects assigned to the control manipulation; thereby the number of subjects in each block respects the overall unequal (and equal) allocation. It is assumed that $n_B$, $ \frac{n_T}{2n} n_B$ and $ \frac{n_C}{2n} n_B$ are integers. The block design with the smallest number of blocks ($B=1$) is equivalent to i/BCRD.

It is clear that $BL(B)$ satisfies Assumptions~\ref{ass:constrained_balance} and \ref{ass:equal_chance_T}. The computation in Appendix \ref{app:ibcrd_minimax} applied to each block implies that $\bSigmaw$ is block diagonal with blocksize $n_B$ where the blocks are all equal to $r\rcontrol (n_B \I_{n_B} - \J_{n_B})/(n_B - 1)$ where $\J_{n_B}$ denotes the $n_B \times n_B$ matrix of all ones.

We show in Theorem~\ref{thm:ibcrd_minimax_blocks} that i/BCRD is minimax among blocking designs, for all response types considered. Note that i/BCRD is a block design with one block ($B=1$, $n_B = 2n$). The proof is found in Section~\ref{app:ibcrd_minimax_blocks} in the Supplemental Material. 

\begin{theorem}\label{thm:ibcrd_minimax_blocks} Consider the responses listed in Table \ref{tab:simple_models} and the definition of $\curlyV_M$ given in Equation~\ref{eq:D_M}, then we have that

\beqn
\mathrm{i/BCRD}  = \argmin_{\W \in BL(B)} \max_{\v \in \curlyV_M} \v^\top \bSigmaw \v.
\eeqn
\end{theorem}

\noindent In summary, if one wishes to gauge design performance by the \emph{worst case MSE}, there is no free lunch; there is no way to mitigate observed covariate imbalances as they will expose your estimator to tail events in the unobserved covariates.

\subsection{The Mean MSE}\label{subsec:mean_criterion}

The minimax criterion forces a design that is extremely conservative to protect against the most adversarial $\v$ which is one arbitrary direction in $\curlyV_M$, a large $2n$-dimensional space. Why not instead design with the \qu{average $\v$} in mind? As the $\x_i$'s are fixed (and thus the $\muvec_T$ and $\muvec_C$ are fixed), we then must consider the $\z_C$ and $\z_T$ as random noise realizations. This is equivalent to assuming that each subject is sampled from a large population also known formally as the \emph{population model} assumption \citep[Chapter 6.2]{Rosenberger2016} whose validity is highly contested especially in clinical trials \citep{Lachin1988}. Regardless of the long-standing debate between the population model and the randomization model, many practitioners assume the population model implicitly and we will likewise in the next two subsections. We reiterate that the mean MSE of this section is one of the most common criterions in the experimental literature.

We formally let $\Z_T$ and $\Z_C$ be the vector random variables that generate the $\z_T$ and $\z_C$ vectors respectively (quantities that are canonically denoted with the Greek letter epsilon). We assume $\expe{\Z_T} = \expe{\Z_C} = \zerovec_{n}$ and that the components in $\Z_T$ are independent, the components of $\Z_C$ are independent and the components of the former and latter are mutually independent. Although the distributional form of both $\Z_T$ and $\Z_C$ are dependent upon the response type specified by $\mathcal{Y}$, we do not make any formal distributional assumptions besides existence of moments.

We now consider the criterion defined as the expectation of the MSE of Equation~\ref{eq:mse} over $\Z_T$ and $\Z_C$. Section~\ref{app:expe_mse} in the Supplementary Material shows this to be

\bneqn\label{eq:expe_mse}
\expesub{\Z_T, \Z_C}{\msesub{\W}{\tauhat}} = \oneover{4n^2}
  \muvec^\top \bSigmaw \muvec + 
  c_Z
\eneqn

\noindent where $\muvec := \muvec_T / r + \muvec_C / \rcontrol$ and $c_Z$ is a constant with respect to the design $\W$ defined in Section~\ref{app:expe_mse}. Intuitively, this expectation operation renders the third term in the MSE decomposition of Equation~\ref{eq:mse_decomposition_three_terms} as a constant and eliminates the second term. Thus, the optimization problem can then be expressed as $\argmin_{\W \in \mathcal{W}} \muvec^\top \bSigmaw \muvec$. Since the second term in $\bSigmaw$ is a constant with respect to the design and the first term can be written as in Equation~\ref{eq:expe_w_w_transpose}, this optimization problem can alternatively be expressed as the design that minimizes 

\bneqn\label{eq:mean_problem}
\muvec^\top \expe{\W\W^\top} \muvec = \sum_{k=1}^{N_{\W}} p_k (\w_k^\top \muvec)^2.
\eneqn


\noindent Thus this design chooses, with equal probability, $\w$'s from the set

\beqn
\argmin_{\braces{\w~:~\w^\top \onevec_{2n} = (r - \rcontrol)n}} (\w^\top \muvec)^2.
\eeqn

In the case of equal allocation, the design that satisfies Assumption~\ref{ass:equal_chance_T} and solves the optimization problem is the design with $N_{\W} = 2$ where the two vectors are $\w_*, -\w_*$ where $\w_*$ is defined to be the unique vector minimizing $\abss{\muvec^\top \w}$, the imbalance between the two arms of the observed covariates' contribution to the responses. This is an NP-hard problem but there are heuristics that could solve it approximately given $\muvec$ as input. In the case of unequal allocation, Equation~\ref{eq:mean_problem} is a difficult problem to solve as the observed covariate imbalances must be minimized while each subject is being placed into the T group in the exact proportion $\frac{n_T}{2n}$ over all $N_{\W}$ vectors. We do not know of a heuristic that can perform this task and we leave its investigation for future work.



However, even in the case of equal allocation, the quantity $\muvec$ is unknown a priori (e.g. because the values of the parameters $\beta_j$ are unknown in the example mean response models given in Table~\ref{tab:simple_models}, column 4). Thus, in real-world scenarios the vector $\w_*$ is absolutely impossible to locate even for the simplest GLMs using heuristic methods. There is one exception --- the limited case where $p = 1$, a linear response model, continuous response and additive treatment effect. Here, the mean response model is 

\beqn
\muvec = \parens{\frac{\beta_0 + \beta_T}{r} + \frac{\beta_0 - \beta_T}{\rcontrol}} \onevec_{2n} + \beta_1 \parens{\oneover{r} + \oneover{\rcontrol}} \x.
\eeqn

\noindent Thus, the optimization problem of Equation~\ref{eq:mean_problem} is solved by the unique vector that minimizes $|\x^\top \w|$. This problem alone can be approximated by numerical methods \citep[Section 2.2.5]{Kapelner2021}. 

In summary, finding the optimal design under the mean MSE criterion is generally intractable. In order to make the problem tractable, we again consider the space of block designs only as in Section \ref{subsubsec:block_designs_worst_case}.


\subsubsection{The Optimal Design Among Block Designs}

Unlike before, we must additionally assume that the order of the elements in the vector $\muvec$. Formally, we assume

\bneqn\label{eq:muvec_order}
\frac{\mu_{T,1}}{r} + \frac{\mu_{C,1}}{\rcontrol} \le \frac{\mu_{T,2}}{r} + \frac{\mu_{C,2}}{\rcontrol} \le \cdots \le \frac{\mu_{T,n}}{r} + \frac{\mu_{C,n}}{\rcontrol}.
\eneqn

\noindent When $p=1$, this assumption is true for all the GLMs given in Table \ref{tab:simple_models} (column 4) since $\mu_{T,i} / r + \mu_{C,i} /\rcontrol$ are monotone in $x_i$ (see also \citet{Kapelner2022b} who discusses this assumption for binary response). When $p>1$, blocks are created that approximate this order. 


Consider the block designs with respect to this order as defined at the end of Section~\ref{subsec:worst_case_criterion}. 
Let $k$ be the smallest possible block size (i.e., $k$ is the smallest integer such that both $\frac{n_T}{2n} k$ and $\frac{n_C}{2n} k$ are integers). Theorem \ref{thm:pm_optimal_blocks} below implies that $BL(k)$ has the smallest mean MSE among block designs. Its proof can be found in Section~\ref{app:opt_expe_mse} of the Supplemental Material.

\begin{theorem}\label{thm:pm_optimal_blocks} Consider the mean criterion given in Equation~\ref{eq:expe_mse}, then

\bneqn\label{eq:blocking_thm_inequality}
\expesub{\Z_T, \Z_C}{\msesub{BL(k)}{\tauhat}} \le \expesub{\Z_T, \Z_C}{\msesub{BL(m)}{\tauhat}},
\eneqn

\noindent for all the vectors $\muvec$ that satisfy the order condition of Equation~\ref{eq:muvec_order} and for any possible block size $m$ where both $\frac{n_T}{2n} m$ and $\frac{n_C}{2n} m$ are integers.
\end{theorem}

\noindent Note that the inequality of Equation~\ref{eq:blocking_thm_inequality} is tight if $2n-1$ of the entries in $\muvec$ are zero and the remaining entry is one, we get equality when $m=1$ corresponding to the the iBCRD design. In the case of equal allocation, the minimum block size is $n_B=k=2$ corresponding to the PM design.

\subsection{The MSE's Tail}\label{subsec:tail_criterion}

The design optimizing for the \emph{worst case MSE} hedges for tail events that have negligible probability and thus is too conservative. The \emph{mean MSE} design (besides being intractable) does not hedge for these tail $\z_t$'s and $\z_C$'s that could introduce large mean squared error with non-negligible probability. Alternatively expressed, the \emph{worst case MSE} minimizes the effect of the third term in the MSE decomposition of Equation~\ref{eq:mse_decomposition_three_terms} while failing to diminish the first term. The \emph{mean MSE} minimizes the first term only because it has the luxury of underhandedly ignoring the third term. We believe these two extant criterions are too crude to measure what we believe to be important: the effect of a large fraction of deleterious tail $\z_t$'s and $\z_C$'s, thereby balancing the contribution of all three terms in the MSE decomposition (Equation~\ref{eq:mse_decomposition_three_terms}).

To properly control this tail MSE, our previous work proposed a novel \emph{tail criterion} which incorporates both the mean and tail events. This new criterion is less conservative than the minimax criterion for more responsive to the tail events. Our previous work showed that the tail criterion leads to a design that is more random than the PB design with two possible allocations $\w_*, -\w_*$ (see end of Section~\ref{subsec:mean_criterion}), and on the other hand, is more balanced and thus more restrictive than i/BCRD. As shown in our previous work and also below, there is a fundamental trade-off: shrinking the mean MSE increases the tail risk of MSE. In this work, we go further and show that block design is asymptotically optimal for the tail criterion under some conditions for many popular response types and unequal treatment allocation. We now define the tail criterion below.

Since MSE is by definition nonnegative, the right tail at a high quantile $q$ (e.g. 95\% or 99\%) can be expressed as

\bneqn\label{eq:actual_quantile_criterion}
\text{Quantile}_{\Z_T, \Z_C}\bracks{\msesub{\W}{\tauhat}, q} = \expesub{\Z_T, \Z_C}{\msesub{\W}{\tauhat}} + c_{\W, q} \sqrt{\varsub{\Z_T, \Z_C}{\msesub{\W}{\tauhat}}}
\eneqn

\noindent where $c_{\W, q}$ is a constant that depends upon the design and the quantile desired. As this constant is not tractable to analyze, we consider the criterion we call the \textit{approximate tail criterion} where we view this constant as fixed with respect to the design,

\bneqn\label{eq:approx_quantile_criterion}
Q_q := \expesub{\Z_T, \Z_C}{\msesub{\W}{\tauhat}} + c_q \sqrt{\varsub{\Z_T, \Z_C}{\msesub{\W}{\tauhat}}}.
\eneqn

\noindent This criterion allows us to incorporate the effect of both the mean and variance with respect to the unobserved covariates and allows us to compare designs asymptotically. We discuss why this approximation is compelling along with alternative criterions that capture both considerations in Section~\ref{sec:discussion}.

Using the expectation from Equation~\ref{eq:expe_mse} and the variance term computed in Section~\ref{app:var_mse} in the Supplemental Material, we can express the approximate tail criterion as 

\bneqn\label{eq:tail_criterion}
Q_q = \oneover{4n^2} \parens{c_Z + \mathcal{B}_1 + c_{q} 
 \sqrt{r^2 \rcontrol^2 \kappa_Z + 4\parens{\mathcal{B}_2 + \mathcal{S}} + 2\mathcal{R} }}
\eneqn

\noindent where the $\kappa_z$ term is a scalar constant with respect to the design $\W$ and a metric similar to the sum of excess kertoses in the components of $\Z_T$ and $\Z_C$ (see Section~\ref{app:var_mse} of the Supplementary Material for its formal definition) and the four other terms are functions of the design described below.\\

\begin{itemize}
\item[$\mathcal{B}_1 :=$] $\muvec^\top \bSigmaw \muvec$ 
        ~~~  
        This is the active term in the MSE expression of Equation~\ref{eq:expe_mse} which was explained as a squared difference or \qu{imbalance} term capturing how different the observed covariates' contribution to the response are between treatment and control groups. This term is unaffected by the distribution of the unobserved $\Z_T$ and $\Z_C$.
        
\item[$\mathcal{B}_2 :=$] $\muvec^\top \bSigmaw \bSigmaZ \bSigmaw \muvec$
        ~~~
        where $\bSigmaZ := \var{\Z_T} / r + \var{\Z_C} / \rcontrol$. This term is a squared difference or imbalance term weighted by the subject-specific variance terms. If all variances of the unobserved effects are bounded, this behaves similarly to term $\mathcal{B}_1$ as it is bounded above by a constant times $\mathcal{B}_1$ (see the proof of Theorem \ref{thm:Q} below).
    
\item[$\mathcal{S} :=$] $r \rcontrol \muvec^\top \bSigmaw \gammaZ$
        ~~~
        where $\gammaZ := \expe{Z_T^3} / r + \expe{Z_C^3} / \rcontrol$, a constant with respect to the design. The term $\mathcal{S}$ measures an average difference or imbalance between the responses in the two groups (a quantity not squared and hence has a sign) and weighted by the subject-specific skewnesses. If the distribution of the $\Z_T$ and $\Z_C$ are approximately symmetric, this term vanishes. Theorem \ref{thm:Q} below shows that $\mathcal{S}$ is asymptotically negligible for block designs. 
    
\item[$\mathcal{R} :=$] $\tr{(\bSigmaw\bSigmaZ)^2}$
        ~~~~
        If all variances of the unobserved effects are the same, this term reduces to a constant times $\normsq{\bSigmaw}_F$ which equals $n$ when $\W$ = CRD (indicating maximal randomness) and equals $n^2$ when $\W$ is PB, the nearly deterministic mirrored pair of vectors (this result is found in Section \ref{subsubsec:tail_PB}). Thus, this term is a measure of the degree of randomness of the design. The $\bSigmaZ$ factor emphasizes or de-emphasizes randomness in certain subjects weighted by the subject-specific variance terms of $\Z_T$ and $\Z_C$. For example, in the case of a binary response, subjects with $\mu_i \approx 50\%$ would have higher values of the variance of their unobserved terms than other subjects resulting in a higher contribution for their degree of randomness in this overall term.\\
\end{itemize}


\noindent It is helpful during our theoretical exploration to define a shifted and scaled quantity of the tail criterion $\tilde{Q}_q := 4n^2 Q - c_Z$ which can be written in the following intuitive expression:

\bneqn\label{eq:Qtilde}
\tilde{Q}_q = \mathcal{B}_1 + c_q \times \sqrt{r^2 \tilde{r}^2 \kappa_Z + 4\parens{\mathcal{B}_2 + \mathcal{S}} + 2 \mathcal{R} }.
\eneqn

\noindent The minimization of $Q_q$ is equivalent to minimization of $\tilde{Q}_q$. However, this minimization

\beqn
\W_*\,|\,q := \argmin_{\W \in \mathcal{W}}\{\tilde{Q}_q\}
\eeqn 

\noindent i.e., the optimal design for the tail criterion as a function of $q$, is \textit{impossible} to compute as the quantities $c_Z, c_q, \kappa_Z, \muvec_T, \muvec_C, \bSigmaZ$ and $\gammaZ$ are unknown. And even if they were to be known, the problem would still be computationally intractable. 

We instead focus on an asymptotic analysis of the terms of $\tilde{Q}_q$. For all designs, $c_Z = \kappa_Z = O(n)$ as they are sums of subject-specific constants. The terms $\mathcal{B}_1$, $\mathcal{B}_2$, $\mathcal{S}$ and $\mathcal{R}$ are design-specific and thus under the experimenter's control. 

Our results are organized as follows. We first locate the asymptotic lower bound of $\tilde{Q}_q$ providing a definition of asymptotic optimality in Section~\ref{subsubsec:lower_bound}. We then asymptotically analyze select designs beginning with i/BCRD (representing \qu{pure} randomness) in Section~\ref{subsubsec:tail_BCRD} and PB (representing \qu{pure} deterministic optimization) in Section~\ref{subsubsec:tail_PB}. These two designs were found to be optimal for the worst case MSE criterion (Section~\ref{subsec:worst_case_criterion}) and the mean MSE criterion (Section~\ref{subsec:mean_criterion}) respectively but we show neither are optimal when gauged by the tail criterion $\tilde{Q}_q$. We then prove that the optimal tail design must lie between these two extremes of pure randomness and pure optimization in Section~\ref{subsubsec:tail_harmony}. We then prove in Section~\ref{subsubsec:tail_block} that blocking design with number of blocks on the order between $n^{1/4}$ and $n$ achieves the lower bound and hence is an asymptotically optimal design.

\subsubsection{An Asymptotic Lower Bound for the Tail Criterion for Any Design}\label{subsubsec:lower_bound}

To prove the lower bound, we first assume the following standard three boundedness conditions of the moments of the $Z_i$'s. Recall that $\gammaZ := \expe{\Z_T^3} / r + \expe{\Z_C^3} / \rcontrol$ and $\bSigmaZ := \var{\Z_T} / r + \var{\Z_C} / \rcontrol$. Let $\brho$ denote the vector of the diagonal entries in $\bSigmaZ$.

\begin{assumption}\label{ass:boundedness1}
The entries of ${\brho}$ are bounded above by  $\bar{C}_{\brho} $ and below by $\underline{C}_{\brho}$ for some positive constants $\underline{C}_{\brho}, \bar{C}_{\brho}$  and  $\frac{1}{2n} \| {\brho}\|^2 \to \bar{\brho} $ for some positive constant $\bar{\brho}$.    
\end{assumption}

\begin{assumption}\label{ass:boundedness2}
$\frac{1}{2n}  \| \gammaZ \|^2 \le C_{\boldsymbol{\gamma}}$ for some positive constant $C_{\boldsymbol{\gamma}}$.    
\end{assumption}

\begin{assumption}\label{ass:boundedness3}
$\frac{1}{2n}  \kappa_Z \to \bar{\kappa}_Z$ for some constant $\bar{\kappa}_Z$.    
\end{assumption}

\noindent In Section~\ref{app:lower} of the Supplementary Material, we prove the following result.

\begin{theorem}\label{thm:lower}
Under Assumptions \ref{ass:boundedness1}--\ref{ass:boundedness3},

\begin{equation}\label{eq:lower}
\liminf_{n \to \infty} \frac{1}{\sqrt{n}} \tilde{Q}_q \ge c_q r \tilde{r}\sqrt{ 2(  \bar{\kappa}_Z +2\bar{\brho})}.
\end{equation}

\end{theorem}

To get some intuition about the lower bound of Equation \ref{eq:lower}, notice that for i/BCRD, $\frac{\mathcal{R}}{n} \to 2 r^2 \tilde{r}^2\bar{\rho}$. Thus, a design attains the lower bound of Equation \ref{eq:lower}, if it is asymptotically as random as i/BCRD, and the balance terms $\mathcal{B}_1$ and $\mathcal{B}_2$ (as well as $\mathcal{S})$ are asymptotically negligible. In Theorem \ref{thm:Q} below we show that block designs with a few number of blocks attains the lower bound under some additional conditions.

\subsubsection{The Tail Criterion under i/BCRD}\label{subsubsec:tail_BCRD}

By the arguments in Section~\ref{app:asymp} of the Supplementary Material, $\mathcal{B}_1 = O(n)$, $\mathcal{R} = O(n)$ (Equation~\ref{eq:Dcomputation}), $\mathcal{B}_2 = O(n)$ (Equation~\ref{eq:Bbound}) and $|\mathcal{S}| = O(n)$ (Equation~\ref{eq:Cbound}) where the first two results are also found in Section 2.2.7 of our previous work for the homoscedastic case. It follows that under i/BCRD the dominant term in $\tilde{Q}$ (Equation~\ref{eq:Qtilde}) is $\mathcal{B}_1$, which is $O(n)$.

\subsubsection{The Tail Criterion under PB}\label{subsubsec:tail_PB}

In the case of equal allocation, the $PB$ design is composed of a single $\w_*$ and its mirror $-\w_*$ that both uniquely minimize $|(\muvec_T + \muvec_C)^\top\w|$. Here, $\mathcal{B}_1 = O(n^3 2^{-2n})$ as shown in \citet[Section 3.3]{Kallus2018} and $\mathcal{R} = O(n^2)$ as shown in Section~\ref{app:harmony} of the Supplementary Material. Furthermore, $\mathcal{B}_2$ is of the same order of $\mathcal{B}_1$ (Equation~\ref{eq:Bbound}) and $|\mathcal{S}| \le \sqrt{O(n) \mathcal{B}_1}$ (Equation~\ref{eq:Cbound}) as also shown in Section~\ref{app:harmony} of the Supplementary Material. Thus, under PB the dominant term in $\tilde{Q}$ (Equation~\ref{eq:Qtilde}) is $\sqrt{\mathcal{R}}$, which is $O(n)$. Thus, PB has the same asymptotic performance as i/BCRD.

\subsubsection{The Asymptotically Optimal Design Harmonizes These Two Extreme Designs}\label{subsubsec:tail_harmony}

In Section~\ref{app:harmony} of the Supplementary Material, we prove the following result.

\begin{theorem}\label{thm:H}
Under Assumptions \ref{ass:boundedness1}--\ref{ass:boundedness3} and under equal allocation, we have that under both i/BCRD and PB, $\frac{1}{\sqrt{n}}\tilde{Q}_q \to \infty$ 
\end{theorem}

Theorem \ref{thm:H} implies that the asymptotically optimal tail criterion design of (Equation~\ref{eq:tail_criterion}), if it exists, must be less imbalanced (and thereby less random) than i/BCRD (i.e. the orders of $\mathcal{B}_1$ and $\mathcal{B}_2$ must increase slower than the rate $n$) and more imbalanced (and thereby more random) than PB (i.e. the order of $\mathcal{R}$ must increase slower than the rate $n^{2}$).

\subsubsection{The Asymptotically Optimal Block Design}\label{subsubsec:tail_block}

We now study the quantile MSE criterion under the set of block designs with number of blocks $B$ similar to the previous two sections. We aim to locate the optimal $B$ as a function of $n$ that minimizes (asymptotically) the quantile criterion.

To prove our results we make two assumptions about blocking designs and provide intuition that motivates these assumptions. 

\begin{assumption}\label{ass:blocking1}
$\mathcal{B}_1 \approx n/B^2$ and if $\frac{B}{n^{1/4}}\to 1$ then $\mathcal{B}_1 \to C_{\mathcal{B}_1}$ for some constant $C_{\mathcal{B}_1}$.
\end{assumption}

\noindent The notation $a_n \approx b_n$ denotes that the limit of $a_n/b_n$ is finite and positive. Assumption \ref{ass:blocking1} connects the number of blocks in the design and the balance term $\mathcal{B}_1$. As shown below, it holds if the vector $\muvec$ is approximately equally spaced, meaning that $\mu_j \approx M/B$ (which is reasonable under the assumption that that $\muvec$ is ordered; see Equation~\ref{eq:muvec_order}). 
Assumption~\ref{ass:blocking1} is motivated by noticing that 

\begin{equation}\label{eq:assump}
\sum_{i=1}^{n_B} (\mu_{b,i} - \bar{\muvec}_b)^2 \approx \frac{n}{B^3}
\end{equation}

\noindent where $\muvec_b$ is the sub-vector of $\muvec$ that corresponds to the $b$-th block is a reasonable assumption. This is because for $b=1$ (and similarly for the other block sizes) $\mu_{1,1},\ldots,\mu_{1,n_B}$ are in $[0, \mu_{n_B}]$ and $\mu_{n_B} \approx M/B$ if the $\mu$'s are ordered approximately uniformly in $[0,M]$.  Therefore, the sample variance satisfies

\[
\frac{1}{n_B}\sum_{i=1}^{n_B} (\mu_{1,i} - \bar{\muvec}_1)^2 \approx \oneover{B^2},
\]

\noindent and since $n_B=2n/B$, Equation~\ref{eq:assump} follows. Equation~\ref{eq:assump} implies Assumption \ref{ass:blocking1} because under block designs we have that 

\[
\mathcal{B}_1=\muvec^\top \bSigmaw \muvec=r \tilde{r} \sum_{b =1}^B \muvec_b^\top \B \muvec_b, 
\]

\noindent where $\B=\frac{n_B}{n_B - 1} \I_{n_B} - \oneover{n_B - 1} \J_{n_B}$. We have that $\B\onevec=\zerovec$, therefore, 

\[
\muvec_b^\top \B \muvec_b = (\muvec_b-\onevec \bar{\muvec}_b)^\top \B (\muvec_b-\onevec \bar{\muvec}_b)=
\frac{n_B}{n_B-1} \sum_{i=1}^{n_B} (\mu_{b,i} - \bar{\muvec}_b)^2 \approx \frac{n_B}{n_B-1} \frac{n}{B^3}, 
\]

\noindent where the last approximation comes from Equation~\ref{eq:assump}. Hence,

\[
\mathcal{B}_1= r \tilde{r} \sum_{b =1}^B \muvec_b^\top \B \muvec_b \approx B \frac{n_B}{n_B-1} \frac{n}{B^3}  \approx \frac{n}{B^2},
\]

\noindent where the last equality is true because $1 \le \frac{n_B}{n_B-1} \le 2$.

\begin{assumption}\label{ass:blocking2}
$\frac{B  \brho^T \bSigmaw \brho}{n^2} \to 0 $. 
\end{assumption}

\noindent When $B$ is of smaller order than $n$, then Assumption~\ref{ass:blocking2} holds true because

\[
\frac{B \brho^T \bSigmaw \brho}{n^2} \le \frac{B  \lambda_{max}(\bSigmaw) \|\brho \|^2 }{n^2} = \frac{B r \tilde{r} \frac{n_B}{n_B-1} \|\muvec\|^2 }{n^2} \le  \frac{2 r \tilde{r} B  \|\muvec\|^2 }{n^2} \to 0.
\]

\noindent Thus, Assumption \ref{ass:blocking2} is only needed to analyze the case where $B/n \to \alpha \in (0,1]$. In this case, the $\mu$'s are balanced by the design in the sense that $\muvec^T \bSigmaw \muvec \approx 1/n$. If the second moments $\brho$ are correlated with the first moments $\muvec$, then it is expected that $\rho^T \bSigmaw \rho$ will be small. In fact it is enough to require that $\rho^T \bSigmaw \rho = o(n)$ for Assumption \ref{ass:blocking2} to hold. In the homoskedastic case, Assumption \ref{ass:blocking2} is always true because then $\brho=\sigma^2 \onevec$ and thus $\brho^T \bSigmaw \brho = 0$. The main result is now given and is proven in Section~\ref{app:asymp} of the Supplementary Material.
 
\begin{theorem}\label{thm:Q}
Under Assumptions \ref{ass:boundedness1}--\ref{ass:blocking2},

\begin{enumerate}[(i)]
    \item If $\frac{B}{n^{1/4}} \to 0 $, then $\frac{1}{\sqrt{n}} \tilde{Q}_q  \to \infty$.
    \item If $\frac{B}{n^{1/4}} \to 1 $, then $\frac{1}{\sqrt{n}} \tilde{Q}_q  \to  C_{\mathcal{B}_1} + c r \tilde{r} \sqrt{ 2(  \bar{\kappa}_Z +2\bar{\brho}) }$.
    \item If $\frac{B}{n^{1/4}} \to \infty$ and $\frac{B}{n} \to 0$, then $\frac{1}{\sqrt{n}} \tilde{Q}_q  \to   c r \tilde{r} \sqrt{ 2(  \bar{\kappa}_Z +2\bar{\brho}) }$.
    \item If $\frac{B}{n} \to \alpha$, for $\alpha \in (0,1]$, then 
    $\frac{1}{\sqrt{n}} \tilde{Q}_q  \to   c r \tilde{r} \sqrt{ 2\left(  \bar{\kappa}_Z +2\bar{\brho} \frac{1}{1-\alpha/2} \right) }$.
\end{enumerate}

\end{theorem}

Theorem \ref{thm:Q}(iii) implies that block designs with $B/n^{1/4} \to \infty$  and $B/n \to 0$ are asymptotically optimal for the tail criterion i.e. it obtains the lower bound of Equation~\ref{eq:lower}.  

Assumption \ref{ass:blocking1} is reasonable when the $\mu_i$'s are ordered as in Equation~\ref{eq:muvec_order} and are approximately uniform. If this is not the case, it is expected that $\mathcal{B}_1$ decreases slower with $B$, and then the optimal block size might be of larger order than $n^{1/4}$.  For example, if  $\mathcal{B}_1 \approx n/B$ (instead of $n/B^2$), then similar arguments used to prove Theorem \ref{thm:Q} will instead imply that the optimal block size $B$ satisfies $B/n^{1/2} \to \infty$ and $B/n \to 0$. More generally, the results above imply if the $\mu_i$'s are neither ordered nor approximately uniform, in order for block designs to perform optimally, it is enough to require that there exists a block design with $B/n \to 0$ such that its balance term $\mathcal{B}_1$ is $o(\sqrt{n})$ as summarized in Remark~\ref{rem:weak} below.

\begin{remark}\label{rem:weak}
Assumptions \ref{ass:blocking1} and \ref{ass:blocking2}, can be relaxed in the sense, that under Assumptions \ref{ass:boundedness1}--\ref{ass:boundedness3}, if there exists a blocking design with $B/n \to 0$ such that its balance term $\mathcal{B}_1$ is $o(\sqrt{n})$, then it asymptotically achieves the lower bound of Equation~\ref{eq:lower} and hence it is asymptotically optimal.
\end{remark}

To interpret the condition that $\mathcal{B}_1=o(\sqrt{n})$, notice that

\[
\frac{1}{2n} \mathcal{B}_1=  r \tilde{r} \frac{1}{B} \sum_{b=1}^B \frac{1}{n_B-1}  \sum_{i=1}^{n_B} (\mu_{b,i} - \bar{\muvec}_b)^2;
\]

\noindent the latter expression is the average variance within blocks multiplied by the constant $r \tilde{r}$. Thus, the condition $\mathcal{B}_1=o(\sqrt{n})$ means that the average variance within blocks is $o(1/\sqrt{n})$. For example, when $B =O(n^{3/4})$ then under Assumption \ref{ass:blocking1}, we have that the average variance within blocks is $O(1/n^{3/2})$ and the requirement in Remark \ref{rem:weak} is weaker, namely that it is $o(1/\sqrt{n})$.



In summary, optimizing for the tail criterion can be accomplished by using a block design with only a few blocks relative to the sample size. 



\section{Simulations}\label{sec:simulations}

\subsection{Setup}

We simulate under a variety of scenarios to verify the performance of the simple estimator (Equation~\ref{eq:estimator}) under the tail criterion (Equation~\ref{eq:tail_criterion}) comports with the results of Theorem~\ref{thm:Q} for blocking designs of varying block sizes across all response types. We simulate under both the equal allocation setting and the unequal allocation setting (specifically, we let $n_C / n_T = 2$). 


The sample size simulated is $2n = 96$ so that we have both a realistically-sized experiment and a wide variety of homogeneously-sized blocking designs, $B \in \braces{1, 2, 3, 4, 6, 8, 12, 16, 24, 48}$ for equal allocation and $B \in \braces{1, 2, 4, 8, 16, 32}$ for the unequal allocation of two controls for every treatment. Note that i/BCRD is equivalent to $B = 1$ and PM is equivalent to $B=n$. We simulate under a number of different covariates $p \in \braces{1, 2, 5}$. Construction of the block designs with $p=1$ are very different from $p>1$ and are explained in detail in Sections~\ref{subsec:block_designs_p_equals_1} and \ref{subsec:block_designs_p_gr_1}. We use all five common response types: continuous, incidence, proportion, count and survival. All response mean models are GLMs and hence contain an internal linear component $\mu_i := \beta_0 + \beta_1 x_{i,1} + \ldots + \beta_p x_{i,p}$ (see columns 3 and 4 of Table~\ref{tab:simple_models}). The values of the covariate coefficients are kept constant for all response types: all $\bbeta = \bracks{+1/5~-1/5~+1/5~-1/5~+1/5}$ and $\beta_0 = -1/5$. The treatment effect is constant at $\beta_T = 1$. One set of fixed covariates are drawn for each pair of response type and number of covariates $p$. Covariate values are always drawn iid and response values are always drawn independently. The covariates were drawn iid uniform. The exact covariate distributions, response distributions and other response parameters are found in Table~\ref{tab:covariate_and_response_settings}. An additional simulation using a long-tailed covariate distribution can be found in Section~\ref{app:additional_simualtions} of the Supplementary Material.  

\begin{table}[ht]
    \small
    \begin{tabular}{llll}
       Response & Covariate & Response & Response\\
       type name &  Distribution &  Distribution & Parameters  \\ \hline
       
       Continuous &  $U(-1, 1)$ & $\normnot{\mu_i}{\sigma^2}$ &          $\sigma = 1$\\ 
       
       Incidence &   $U(-3, 3)$ & $\bernoulli{\mu_i}$ &                   N/A  \\
       
       Proportion &  $U(-1, 1)$ & $\betanot{\phi\mu_i}{\phi(1-\mu_i)}$ & $\phi = 2$\\ 
       
       Count &       $U(-1, 1)$ & $\poisson{\mu_i}$ & N/A \\
       
       Survival &    $U(-1, 1)$ & $\weibullnot{\mu_i / \Gamma\parens{1 + 1/k}}{k}$ 
& $k = 4$\\ \hline
    \end{tabular}
    \caption{Simulation settings and parameters.}
    \label{tab:covariate_and_response_settings}
\end{table}



The tail criterion is computed on the distribution of the MSE over $\Z_T, \Z_C$. To approximate this quantile, we generate $N_y := 100,000$ pairs $\y_T, \y_C$ using draws from the response distribution (and hence implicitly draw pairs $\z_T, \z_C$). Using these draws, we compute $\v$ and then can exactly compute the MSE via the quadratic form of Equation~\ref{eq:mse}. (The entries of $\bSigmaw$ are for each block design specified by the specific $B$). We approximate the $q=.95$ tail criterion using the empirical 95th percentile of the $N_y$ MSE values. And we also calculate the approximate tail when $c_q = 1.645$, the 95\%ile of the standard normal in order to compare the actual quantile (Equation~\ref{eq:actual_quantile_criterion}) to our approximate tail criterion (Equation~\ref{eq:approx_quantile_criterion}).



\subsection{Block design construction for $p = 1$}\label{subsec:block_designs_p_equals_1}

Simply ordering the subjects by the order of the one covariate measurement then splicing into $B$ even sets provides the block design structure as all observed covariate response functions are monotonic in $\x$ (see columns 3 and 4 of Table~\ref{tab:simple_models}). Thus, the entries of $\muvec$ are ordered within the blocks corresponding to the block-diagonal structure of $\bSigmaw$. 
Simulation results for $p=1$ constitute empirical confirmation of our mathematical results.

\subsection{Block design construction for $p > 1$}\label{subsec:block_designs_p_gr_1}

For $p>1$, we cannot order $\muvec$ by using the covariate data (as $\bbeta$ is considered unknown in the real world for any response type). Thus, we employ blocking as done in practice. We arbitrarily consider the first two covariates (we do not choose to use more covariates as the number of blocks increases exponentially in $p$, a common problem in experimental practice). We first order the subjects by the first covariate. Then within blocks of size $2 n_B$, we sort by the second covariate. Thus, we emphasize, simulation results for $p > 1$ serve as an empirical check on the robustness of our theoretical results as $\muvec$ will not be exactly sorted as required by Equation~\ref{eq:muvec_order}.

\subsection{Results}

\begin{figure}[htp]
\centering
     \begin{subfigure}[b]{\textwidth}
         \centering
         \includegraphics[width=6in]{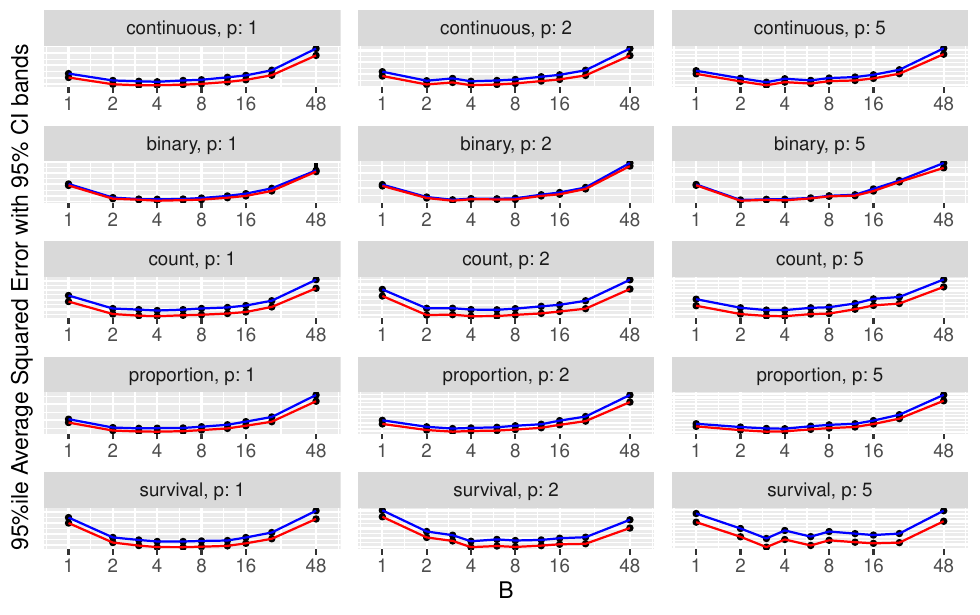}
         \caption{Equal Allocation}
         \label{fig:all_response_q_95_equal_allocation}
     \end{subfigure}\\
     \begin{subfigure}[b]{\textwidth}
         \centering
         \includegraphics[width=6in]{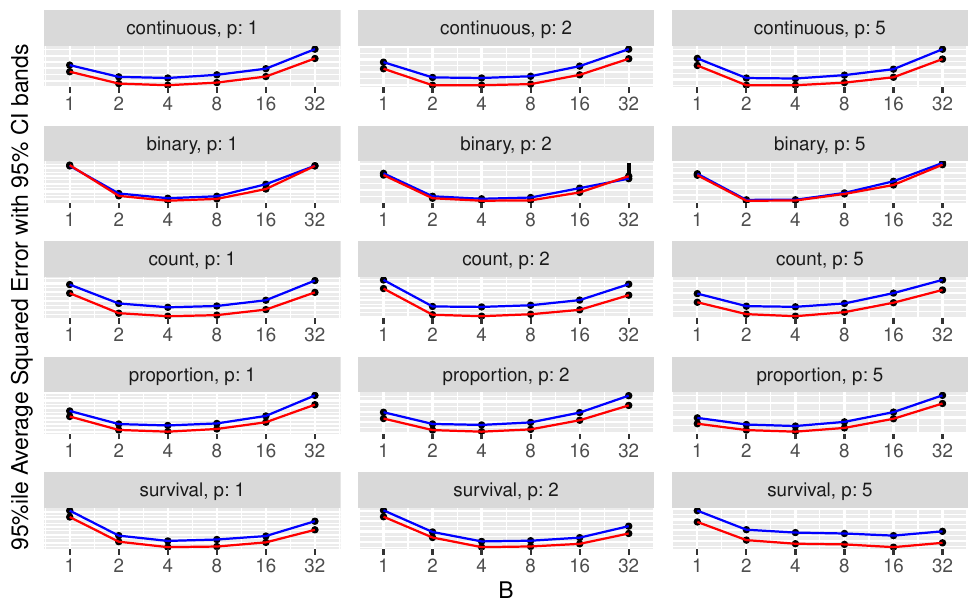}
         \caption{Unequal Allocation ($n_C:n_T$ = 2:1)}
         \label{fig:all_response_q_95_unequal_allocation}
     \end{subfigure} 
\caption{Simulation results for mean-centered uniform covariate realizations, $n=96$, $p = 1,2,5$, all response types and all appropriate number of blocks. The y-axis is relative to simulation and thus its values are unshown. Blue lines are empirical 95\%ile MSE Tails of (Equation~\ref{eq:actual_quantile_criterion}) and red lines are approximate quantiles computed via the average MSE plus 1.645 times the empirical standard error of the MSE (Equation~\ref{eq:approx_quantile_criterion}).}
\end{figure}

The $q=.95$ results for equal allocation are presented in Figure~\ref{fig:all_response_q_95_equal_allocation}. We can see that for $p = 1$ and all response types, $B \in \braces{2, 4, 8}$ are the best-performing designs. As $n^{1/4} \approx 3.4$ and $n^{1/2} \approx 11.3$, four blocks and eight blocks reasonably comport with this order of magnitude. This empirical result also vindicates the assumptions made about the blocking design (Assumptions~\ref{ass:blocking1} and \ref{ass:blocking2}) and additionally demonstrates the asymptotic results apply in the realistic experimental sample size of $n \approx 100$. 

For $p > 1$, 
the $\mu_i$'s are not precisely ordered into blocks, but only approximately by the arrangement of blocks of two covariates. Notwithstanding, the results for $p > 1$ overall demonstrate the robustness of our theoretical results to situations of \qu{imperfect} blocking; and imperfect blocking is the design employed in all real-world experimental settings. Thus, the $p > 1$ results are an empirical demonstration of Remark~\ref{rem:weak} which provides evidence that our theory applies to real-world block designs that imperfectly order $\muvec$.

Figure~\ref{fig:all_response_q_95_unequal_allocation} shows the results for unequal allocation where treatment is administered to $1/3$ of the subjects. Here, we see the same picture as in the equal allocation simulation for both $p=1$ and $p > 1$ further supporting our theoretical results in real-world scenarios.

For both equal and unequal allocation, the approximate tail criterion computed using the sample average, sample standard deviation and $c_q = 1.645$ (the standard normal 95\% quantile) is nearly the same as the empirical 95\%ile (compare red and blue lines for all plot windows in Figures~\ref{fig:all_response_q_95_equal_allocation} and \ref{fig:all_response_q_95_unequal_allocation}). This comparison demonstrates that the approximate tail criterion (Equation~\ref{eq:approx_quantile_criterion}) is nearly equal to the actual quantile tail criterion (Equation~\ref{eq:actual_quantile_criterion}) for the gamut of real world scenarios in all response types. Whatever difference exists between the approximate quantile and the actual quantile further does not affect the optimal lowest MSE designs (visualized as the $B$ value corresponding to the minimum of the curves).

\section{Discussion}\label{sec:discussion}

Randomized experimentation comes with a choice to the experimenter: how to randomize exactly? Should this choice differ depending on the experimental response metric? One extreme position is full randomization independent of subjects' covariate imbalance across arms. Here, there can be high observed covariate imbalance among the arms translating to a highly visible risk of bias. However, we prove herein that if the experimenter is concerned with the worst possible configuration of unobserved covariates, complete randomization is the optimal design for all common experimental response types. Another extreme is to find one perfect assignment that minimizes covariate imbalance across arms. However, we prove herein that \emph{only} if one assumes we can average over the unmeasured covariates, this deterministic allocation would be the optimal design for all response types under equal allocation. But this is still unrealistic in practice; that optimal allocation would be impossible to locate as for the vast majority of covariate-response models, knowing the covariate parameters and functional form is required. Even if the parameters and functional form were to be known, the problem would be NP-hard and thus can only be approximated.

Other commonly employed designs are in the \qu{middle} of these two extremes; they sacrifice randomness for imbalance reduction. How to harmonize these two considerations to rigorously choose a design was the subject of our previous work which exclusively considered the continuous response and equal treatment allocation setting. Therein, we proposed the \qu{quantile tail criterion}, a metric on designs that reveals how the degree of randomness and degree of imbalance reduction trade off with one another within a design. 

Herein, we derive an \emph{approximate quantile tail criterion} for general response types. We then prove that the asymptotically optimal design is achieved with a block design with very few blocks --- on the order of between $n^{1/4}$ and $n$ as seen in case (iii) of Theorem \ref{thm:Q}. Blocking with few blocks is an intuitive imbalance-reducing idea proposed by Fisher nearly 100 years ago and is one of the most popular designs employed today in clinical trials. This work vindicates this common practice for all common response types by analyzing its performance from a new angle. We  demonstrate that blocking is the optimal design harmonizing the degree of randomness and imbalance reduction. Although this work focuses on estimation, these blocking designs should have high power as well as many allocations can be generated which are highly independent \citep{Krieger2020b}.

Our criterion differs from the exact MSE quantile by assuming $c_{\W, q} = c_q$, a constant with respect to the $\W, \Z_T, \Z_C$. However, what if for some pathological designs or unobserved covariate distributions it is not constant? If one has such a concern, our approximate tail criterion can always be reinterpreted as controlling the upper bound on the $q$th quantile as $c_q$ can be set to $\sqrt{q/(1-q)}$, a bound based on Chebyshev's Inequality \cite[Section 3.1]{Bagui2004}. For example, if $q=.95$, letting $c_q = 4.359$ would be a crude upper bound for the exact 95\%ile of MSE. As our results our asymptotic, the value of $c_q$ (no matter how large) does not matter (as it is constant with respect to $\W, \Z_T, \Z_C$ and sample size). In this perspective, we are minimizing how large a destructive event \emph{can potentially be}. And we believe this is still a compelling criterion to control.

Regardless, we believe such pathological settings will be rare as theoretically discussed in \cite{Sharpe1970} and as empirically demonstrated by the simulations of Section~\ref{sec:simulations}. Here, we found that employing $c_q$ as the standard normal quantile provides near-exact quantiles across the block designs considered. We conjecture there is a central limit theorem on $\msesub{\W}{\tauhat}$ over realizations of $\Z_T, \Z_C$ if $N_{\W}$ is large; this can be left to future work.

In Section~\ref{subsec:mean_criterion} we assumed knowledge of the order $\muvec$, which is reasonable when $p=1$ but not generally. When this assumption holds and $\muvec$ is approximately equally spaced (see Assumption~\ref{ass:blocking1}) then Theorem \ref{thm:Q} implies that the optimal number of blocks is between $n^{1/4}$ and $n$. More generally, Remark~\ref{rem:weak} shows that blocking designs with number of blocks of smaller order than $n$ and that are balanced enough are asymptotically optimal (in the sense that they asymptotically achieve the lower bound of Theorem~\ref{thm:lower}). The results of Section~\ref{sec:simulations} demonstrate that in the models we considered, indeed blocking design with only a few number of blocks (e.g., 4 or 8) are optimal (among blocking designs) also for the case $p>1$.  

There are many other avenues to explore. First, although we proved that blocking with few blocks are asymptotically optimal, we did not prove it is uniquely optimal. Other designs may also satisfy the lower bound of Equation~\ref{eq:lower} and offer better finite sample properties. Other optimal designs are additionally of profound importance because blocking suffers from two main practical problems (1) the curse of dimensionality as it can only handle $p \leq 3$ in practice and (2) the design requires arbitrary cut-points if the covariate is continuous. 
Another area of exploration could be understanding the tail criterion for other designs such as rerandomization and greedy pair switching. Additionally, survival is an important experimental endpoint in clinical trials but this work only considered uncensored survival measurement while the vast majority of survivals collected in real-world experiments are censored. Also, we did not investigate the distribution of the $\mu_i$'s relative to the optimal block size and we believe this can be a source of fruitful future work. In our previous work (Section 2.3), we explored the ordinary least squares estimator. This avenue would be difficult to explore for general response types as the analogous estimators do not have closed form (e.g. if the response is binary, the multivariate logistic regression estimates are the result of computational iteration). It is our intuition the theorems found herein would extend to these commonly used covariate-controlled estimators, but more work is necessary. Further, our tail criterion is only one metric that combines mean and standard deviation of MSE. Other criterions can be explored e.g. a tail simultaneous in allocation and unmeasured covariates, the ratio of mean to standard deviation (the coefficient of variation), etc. Also, work can be done to find the optimal imbalance design under unequal allocation, the analogue of PB in equal allocation (see discussion below Equation~\ref{eq:mean_problem}). Finally, in the case of unequal allocation, we can develop a PM analogue by performing $(n_C / n_T) : 1$ matching. If $n_C / n_T$ is not a natural number, some matches will have $\lfloor n_C / n_T \rfloor$ control subjects for one treatment subject and some matches can have $\lceil n_C / n_T \rceil$ control subjects for one treatment subject (in approximately the correct proportion). An algorithm to create these tuplets based on intra-subject observed covariate distance is only known to be able to be solved in polynomial time when each matched set has one treated and one control. Clearly, this more general $(n_C / n_T) : 1$ matching problem is more complex and might very well be NP-hard.

\subsection*{Funding information}

This research was supported by Grant No 2018112 from the United States-Israel Binational Science Foundation (BSF).

\subsection*{Author contributions}

All authors have accepted responsibility for the entire content of this manuscript and approved its submission.

\subsection*{Conflict of interest}

Authors state no conflict of interest.

\subsection*{Data availability statement}

All data generated or analysed during this study are included in this published article and its supplementary information files.


\begin{thebibliography}{}

\bibitem[Bagui and Bhaumik, 2004]{Bagui2004}
Bagui, S.~C. and Bhaumik, D.~K. (2004).
\newblock Glimpses of inequalities in probability and statistics.
\newblock {\em International Journal of Statistical Sciences}, 3:9--15.

\bibitem[Bertsimas et~al., 2015]{Bertsimas2015}
Bertsimas, D., Johnson, M., and Kallus, N. (2015).
\newblock The power of optimization over randomization in designing experiments involving small samples.
\newblock {\em Operations Research}, 63(4):868--876.

\bibitem[Chattopadhyay et~al., 2022]{Chattopadhyay2022}
Chattopadhyay, A., Morris, C.~N., and Zubizarreta, J.~R. (2022).
\newblock Balanced and robust randomized treatment assignments: the finite selection model.
\newblock {\em arXiv preprint arXiv:2205.09736}.

\bibitem[Cornfield, 1959]{Cornfield1959}
Cornfield, J. (1959).
\newblock Principles of research.
\newblock {\em American journal of mental deficiency}, 64:240--252.

\bibitem[Dumville et~al., 2006]{Dumville2006}
Dumville, J., Hahn, S., Miles, J., and Torgerson, D. (2006).
\newblock The use of unequal randomisation ratios in clinical trials: a review.
\newblock {\em Contemporary clinical trials}, 27(1):1--12.

\bibitem[Efron, 1971]{Efron1971}
Efron, B. (1971).
\newblock {Forcing a sequential experiment to be balanced}.
\newblock {\em Biometrika}, 58(3):403--417.

\bibitem[Feng and Hu, 2018]{Feng2018}
Feng, C. and Hu, F. (2018).
\newblock Optimal responses-adaptive designs based on efficiency, ethic, and cost.
\newblock {\em Statistics and Its Interface}, 11(1):99--107.

\bibitem[Fisher, 1925]{Fisher1925}
Fisher, R.~A. (1925).
\newblock {\em Statistical methods for research workers}.
\newblock Edinburgh Oliver \& Boyd.

\bibitem[Fisher, 1935]{Fisher1935}
Fisher, R.~A. (1935).
\newblock {\em The Design of Experiments}.
\newblock Oliver and Boyd, 1 edition.

\bibitem[Freedman, 2008]{Freedman2008}
Freedman, D.~A. (2008).
\newblock {On regression adjustments to experimental data}.
\newblock {\em Advances in Applied Mathematics}, 40(2):180--193.

\bibitem[Greevy et~al., 2004]{Greevy2004}
Greevy, R., Lu, B., Silber, J.~H., and Rosenbaum, P. (2004).
\newblock Optimal multivariate matching before randomization.
\newblock {\em Biostatistics}, 5(2):263--275.

\bibitem[Harshaw et~al., 2019]{Harshaw2019}
Harshaw, C., S{\"a}vje, F., Spielman, D., and Zhang, P. (2019).
\newblock Balancing covariates in randomized experiments with the gram--schmidt walk design.
\newblock {\em arXiv preprint arXiv:1911.03071}.

\bibitem[Harville, 1975]{Harville1975}
Harville, D.~A. (1975).
\newblock Experimental randomization: {W}ho needs it?
\newblock {\em The American Statistician}, 29(1):27--31.

\bibitem[Imbens and Rubin, 2015]{Imbens2015}
Imbens, G.~W. and Rubin, D.~B. (2015).
\newblock {\em Causal inference in statistics, social, and biomedical sciences}.
\newblock Cambridge University Press.

\bibitem[Kallus, 2018]{Kallus2018}
Kallus, N. (2018).
\newblock Optimal a priori balance in the design of controlled experiments.
\newblock {\em Journal of the Royal Statistical Society: Series B (Statistical Methodology)}, 80(1):85--112.

\bibitem[Kapelner et~al., 2022]{Kapelner2022b}
Kapelner, A., Krieger, A.~M., and Azriel, D. (2022).
\newblock The role of pairwise matching in experimental design for an incidence outcome.
\newblock {\em arXiv preprint arXiv:2209.00490}.

\bibitem[Kapelner et~al., 2021]{Kapelner2021}
Kapelner, A., Krieger, A.~M., Sklar, M., Shalit, U., and Azriel, D. (2021).
\newblock Harmonizing optimized designs with classic randomization in experiments.
\newblock {\em The American Statistician}, 75(2):195--206.

\bibitem[Kiefer and Wolfowitz, 1959]{Kiefer1959}
Kiefer, J. and Wolfowitz, J. (1959).
\newblock Optimum designs in regression problems.
\newblock {\em The Annals of Mathematical Statistics}, pages 271--294.

\bibitem[Krieger et~al., 2019]{Krieger2019}
Krieger, A.~M., Azriel, D., and Kapelner, A. (2019).
\newblock Nearly random designs with greatly improved balance.
\newblock {\em Biometrika}, 106(3):695--701.

\bibitem[Krieger et~al., 2020]{Krieger2020b}
Krieger, A.~M., Azriel, D., Sklar, M., and Kapelner, A. (2020).
\newblock Improving the power of the randomization test.
\newblock {\em arXiv preprint arXiv:2008.05980}.

\bibitem[Lachin, 1988]{Lachin1988}
Lachin, J.~M. (1988).
\newblock Properties of simple randomization in clinical trials.
\newblock {\em Controlled clinical trials}, 9(4):312--326.

\bibitem[Land and Doig, 1960]{Land1960}
Land, A.~H. and Doig, A.~G. (1960).
\newblock An automatic method of solving discrete programming problems.
\newblock {\em Econometrica: Journal of the Econometric Society}, pages 497--520.

\bibitem[Li et~al., 2018]{Li2018}
Li, X., Ding, P., and Rubin, D.~B. (2018).
\newblock Asymptotic theory of rerandomization in treatment--control experiments.
\newblock {\em Proceedings of the National Academy of Sciences}, 115(37):9157--9162.

\bibitem[Lin, 2013]{Lin2013}
Lin, W. (2013).
\newblock Agnostic notes on regression adjustments to experimental data: Reexamining freedman's critique.
\newblock {\em The Annals of Applied Statistics}, 7(1):295--318.

\bibitem[Morgan and Rubin, 2012]{Morgan2012}
Morgan, K.~L. and Rubin, D.~B. (2012).
\newblock Rerandomization to improve covariate balance in experiments.
\newblock {\em The Annals of Statistics}, pages 1263--1282.

\bibitem[Nam, 1973]{Nam1973}
Nam, J.-M. (1973).
\newblock Optimum sample sizes for the comparison of the control and treatment.
\newblock {\em Biometrics}, pages 101--108.

\bibitem[Nordin and Schultzberg, 2022]{Nordin2022}
Nordin, M. and Schultzberg, M. (2022).
\newblock Properties of restricted randomization with implications for experimental design.
\newblock {\em Journal of Causal Inference}, 10(1):227--245.

\bibitem[Rosenberger and Lachin, 2016]{Rosenberger2016}
Rosenberger, W.~F. and Lachin, J.~M. (2016).
\newblock {\em Randomization in clinical trials: theory and practice}.
\newblock John Wiley \& Sons, second edition.

\bibitem[Rubin, 2005]{Rubin2005}
Rubin, D.~B. (2005).
\newblock Causal inference using potential outcomes: Design, modeling, decisions.
\newblock {\em Journal of the American Statistical Association}, 100(469):322--331.

\bibitem[Ryeznik and Sverdlov, 2018]{Ryeznik2018}
Ryeznik, Y. and Sverdlov, O. (2018).
\newblock A comparative study of restricted randomization procedures for multiarm trials with equal or unequal treatment allocation ratios.
\newblock {\em Statistics in Medicine}, 37(21):3056--3077.

\bibitem[Senn, 2013]{Senn2013}
Senn, S. (2013).
\newblock Seven myths of randomisation in clinical trials.
\newblock {\em Statistics in Medicine}, 32(9):1439--1450.

\bibitem[Sharpe, 1970]{Sharpe1970}
Sharpe, K. (1970).
\newblock Robustness of normal tolerance intervals.
\newblock {\em Biometrika}, 57(1):71--78.

\bibitem[Smith, 1918]{Smith1918}
Smith, K. (1918).
\newblock On the standard deviations of adjusted and interpolated values of an observed polynomial function and its constants and the guidance they give towards a proper choice of the distribution of observations.
\newblock {\em Biometrika}, 12(1/2):1--85.

\bibitem[Steinberg and Hunter, 1984]{Steinberg1984}
Steinberg, D.~M. and Hunter, W.~G. (1984).
\newblock Experimental design: review and comment.
\newblock {\em Technometrics}, 26(2):71--97.

\bibitem[Student, 1938]{Student1938}
Student (1938).
\newblock Comparison between balanced and random arrangements of field plots.
\newblock {\em Biometrika}, pages 363--378.

\bibitem[Sverdlov and Ryeznik, 2019]{Sverdlov2019}
Sverdlov, O. and Ryeznik, Y. (2019).
\newblock Implementing unequal randomization in clinical trials with heterogeneous treatment costs.
\newblock {\em Statistics in Medicine}, 38(16):2905--2927.

\bibitem[Teugels, 1990]{Teugels1990}
Teugels, J.~L. (1990).
\newblock Some representations of the multivariate bernoulli and binomial distributions.
\newblock {\em Journal of multivariate analysis}, 32(2):256--268.

\bibitem[Torgerson and Campbell, 1997]{Torgerson1997}
Torgerson, D. and Campbell, M. (1997).
\newblock Unequal randomisation can improve the economic efficiency of clinical trials.
\newblock {\em Journal of Health Services Research \& Policy}, 2(2):81--85.

\bibitem[Tymofyeyev et~al., 2007]{Tymofyeyev2007}
Tymofyeyev, Y., Rosenberger, W.~F., and Hu, F. (2007).
\newblock Implementing optimal allocation in sequential binary response experiments.
\newblock {\em Journal of the American Statistical Association}, 102(477):224--234.

\bibitem[Wong and Zhu, 2008]{Wong2008}
Wong, W.~K. and Zhu, W. (2008).
\newblock Optimum treatment allocation rules under a variance heterogeneity model.
\newblock {\em Statistics in Medicine}, 27(22):4581--4595.

\bibitem[Wu, 1981]{Wu1981}
Wu, C.~F. (1981).
\newblock On the robustness and efficiency of some randomized designs.
\newblock {\em The Annals of Statistics}, pages 1168--1177.

\bibitem[Zhu and Liu, 2021]{Zhu2021}
Zhu, K. and Liu, H. (2021).
\newblock Pair-switching rerandomization.
\newblock {\em arXiv preprint arXiv:2103.13051}.

\end{thebibliography}

\appendix
\pagebreak

\begin{center}
   \Large{Supplementary Information for} \\ 
   \LARGE{\qu{\ourtitle}} 
\end{center}

\section{Technical Proofs}

\subsection{Estimator Computation}\label{app:estimator_computation}

To compute the estimator we begin with:

\beqn
\tauhat &:=& \YbarT - \YbarC \\
&=& \oneover{n_T} \overtwo{(\onevec_{2n} + \W)^\top} \Y - \oneover{n_C} \overtwo{(\onevec_{2n} - \W)^\top} \Y \\
&=& \half \parens{\frac{(\onevec_{2n} + \W)}{n_T} - \frac{(\onevec_{2n} - \W)}{n_C}}^\top \Y
\eeqn

\noindent We now substitute for $\Y$ from Equation~\ref{eq:potential_outcomes}. To simplify the notation, let $\v := \y_T + \y_C$ which is slightly different than the notation in the main text, $v_0 := \onevec_{2n} \v$, $y_{T,0} := \onevec_{2n}^\top \y_T$ and $y_{C,0} := \onevec_{2n}^\top \y_C$.

\beqn
\tauhat &=& \oneover{4} \parens{\frac{(\onevec_{2n}^\top + \W)^\top}{n_T} - \frac{(\onevec_{2n}^\top - \W^\top)}{n_C}} \parens{\v + \diag{\W}(\y_T - \y_C)} \\
&=& \oneover{4} \parens{
\frac{(v_0+ \W^\top\v)}{n_T} - \frac{(v_0 - \W^\top\v)}{n_C} + \parens{\frac{\onevec_{2n}^\top + \W^\top}{n_T} - \frac{\onevec_{2n}^\top - \W^\top}{n_C}} \diag{\W}(\y_T - \y_C)
} \\
&=& \oneover{4} \parens{
\frac{(v_0+ \W^\top\v)}{n_T} - \frac{(v_0 - \W^\top\v)}{n_C} + \parens{\frac{\W^\top + \onevec_{2n}^\top}{n_T} - \frac{\W^\top - \onevec_{2n}^\top}{n_C}} (\y_T - \y_C)
} \\
&=& \footnotesize\oneover{4} \parens{
\frac{(v_0+ \W^\top\v)}{n_T} - \frac{(v_0 - \W^\top\v)}{n_C} + \frac{\W^\top\y_T +  y_{T,0}}{n_T} - \frac{\W^\top\y_T - y_{T,0}}{n_C} - \frac{\W^\top\y_C + y_{C,0}}{n_T} + \frac{\W^\top\y_C - y_{C,0}}{n_C}
} \\
&=& \footnotesize\oneover{4} \parens{ 
    \frac{
        v_0 + \W^\top\v + \W^\top\y_T + y_{T,0} - \W^\top\y_C - y_{C,0}
    }{n_T} + 
    \frac{
        -v_0 + \W^\top\v - \W^\top\y_T + y_{T,0} + \W^\top\y_C - y_{C,0}
    }{n_C}
} \\
&=& \oneover{4} \parens{ 
    \frac{
        v_0 + y_{T,0} - y_{C,0}
    }{n_T} + 
    \frac{
        -v_0 + y_{T,0} - y_{C,0}
    }{n_C} +
    \frac{
        \W^\top (\v + \y_T - \y_C)
    }{n_T} + 
    \frac{
        \W^\top (\v - \y_T + \y_C)
    }{n_C}
} \\
&=& \oneover{4} \parens{ 
    \frac{
        2y_{T,0}
    }{n_T} + 
    \frac{
        -2y_{C,0}
    }{n_C} +
    \W^\top \parens{
        \frac{
            2\y_T
        }{n_T} + 
        \frac{
            2\y_C
        }{n_C}
    }
} \\
&=& \half \parens{ 
    \frac{
        y_{T,0}
    }{n_T} - 
    \frac{
        y_{C,0}
    }{n_C} +
    \W^\top \parens{
        \frac{
            \y_T
        }{n_T} + 
        \frac{
            \y_C
        }{n_C}
    }
}  \\
&=& \half \parens{ 
    \onevec_{2n}^\top \parens{\frac{
            \y_T
        }{n_T} - 
        \frac{
            \y_C
        }{n_C}
    }+
    \W^\top \parens{
        \frac{
            \y_T
        }{n_T} + 
        \frac{
            \y_C
        }{n_C}
    }
}  \\
&=& \oneover{2n} \parens{ 
    \onevec_{2n}^\top \parens{\frac{
            \y_T
        }{r} - 
        \frac{
            \y_C
        }{\rcontrol}
    }+
    \W^\top \parens{
        \frac{
            \y_T
        }{r} + 
        \frac{
            \y_C
        }{\rcontrol}
    }
}
\eeqn

\subsection{Unbiasedness of the Estimator}\label{app:estimator_unbiasedness}

Under Assumption~\ref{ass:equal_chance_T},

\beqn
\expesub{\W}{\tauhat} &=& \half\parens{ 
    \onevec_{2n}^\top \parens{\frac{
            \y_T
        }{n_T} - 
        \frac{
            \y_C
        }{n_C}
    }+
    \frac{n_T - n_C}{2n}\onevec_{2n}^\top \parens{
        \frac{
            \y_T
        }{n_T} + 
        \frac{
            \y_C
        }{n_C}
    }
} \\
&=& \half\onevec_{2n}^\top \parens{ 
    \frac{
            \y_T
        }{n_T} - 
        \frac{
            \y_C
        }{n_C}
    +
    \frac{n_T - n_C}{2n}\parens{
        \frac{
            \y_T
        }{n_T} + 
        \frac{
            \y_C
        }{n_C}
    }
} \\
&=& \half\onevec_{2n}^\top \parens{ 
    \y_T \parens{\oneover{n_T} + \frac{n_T - n_C}{2n n_T}} +
    \y_C \parens{-\oneover{n_C} + \frac{n_T - n_C}{2n n_C}}
} \\
&=& \half\onevec_{2n}^\top \parens{ 
    \y_T \frac{n + n_T - n_C}{2n n_T} +
    \y_C \frac{-n + n_T - n_C}{2n n_C}
} \\
&=& \half\onevec_{2n}^\top \parens{ 
    \y_T \frac{2n_T}{2n n_T} +
    \y_C \frac{-2n_C}{2n n_C}
} \\
&=& \oneover{2n} \onevec_{2n}^\top \parens{\y_T - \y_C} = \tau~~\blacksquare
\eeqn

\subsection{i/BCRD is Minimax for Continuous Response}\label{app:ibcrd_minimax}

The maximum of the quadratic form of Equation~\ref{eq:worst_case_problem} is $\lambda_{\text{max}}(\bSigmaw) \normsq{\v}$. The form of $\bSigmaw$ for the i/BCRD design is crucial to understand thus we compute it now. 

For two subjects $i$ and $j$ the joint PMF of their assignments is:

\beqn
\prob{W_i = +1, W_j = +1} &=& \frac{n_T(n_T-1)}{2n(2n-1)}\\
\prob{W_i = +1, W_j = -1} = \prob{W_i = -1, W_j = +1} &=& \frac{n_T n_C}{2n(2n-1)}\\
\prob{W_i = -1, W_j = -1} &=& \frac{n_C(n_C-1)}{2n(2n-1)} 
\eeqn

\noindent Letting $b := (n_T - n_C)^2 / (2n)^2$ we first compute the diagonal entries of $\bSigmaw$:

\bneqn\label{eq:diagonal_of_bSigmaW}
\expe{W_i} &=& \frac{n_T - n_C}{2n} = \frac{r - \rcontrol}{2} = \sqrt{b} \nonumber \\
\expe{W_i^2} &=& \frac{n_T}{2n} + \frac{n_C}{2n} = \frac{r + \rcontrol}{2} = 1 \nonumber \\
\var{W_i} &=& 1 - b = \frac{n_T n_C}{4n^2} = r\rcontrol
\eneqn

\noindent Note that the above is valid for any design that satisfies Assumption~\ref{ass:equal_chance_T}. We then compute the off-diagonal elements of $\bSigmaw$ for i/BCRD:

\beqn
\expe{W_i W_j} &=& \frac{n_C(n_C-1) + n_T(n_T-1) - 2 n_T n_C}{2n(2n-1)} = \frac{2nb - 1}{2n-1} \\
\cov{W_i}{W_j} &=& \frac{2nb - 1}{2n-1} - b = \frac{b - 1}{2n - 1}.
\eeqn

\noindent Thus,

\beqn
\bSigmaw = -(b - 1) \parens{\frac{2n}{2n-1} \I_{2n} - \frac{1}{2n-1}\J_{2n}}
= \frac{r \tilde{r}}{2n-1} \parens{2n \I_{2n} - \J_{2n}}.
\eeqn

\noindent Note that when $n_T = n_C$, then $b = 0$ and the matrix reduces to the variance covariance matrix for the BCRD design. Since the BCRD variance covariance has one zero eigenvalue and $2n-1$ eigenvalues with value $\frac{2}{2n-1}$, this matrix being a scalar multiple will have one zero eigenvalue and $2n-1$ eigenvalues with value $\frac{2n(b-1)}{2n-1}$.

Since the sum of the eigenvalues is equal to $\tr{\bSigmaw}$ and since all designs we consider have the same diagonal and $2n-1$ positive eigenvalues, any other design besides i/BCRD will have different eigenvalues hence a larger maximum eigenvalue. Thus the design that minimizes $\lambda_{\text{max}}(\bSigmaw) \normsq{\v}$ is i/BCRD. $\blacksquare$

\subsection{i/BCRD is Minimax for all Response Types}\label{app:ibcrd_minimax_blocks}

For the continuous response case, the result follows from Theorem \ref{thm:iBCRD_minimax_for_continuous}. We now consider the non-continuous response case. Notice that in this case ${\curlyV}_M$ has the form 
\[
\{v \in \curlyV : \min(v_1,\ldots,v_{2n}) \ge 0 ~,~ \max(v_1,\ldots,v_{2n}) \le M \}
\]
(for the incidence and proportion responses, $M=1/r+1/\tilde{r}$). 
The key idea is to show that for block designs

\bneqn\label{eq:key}
 \max_{\v \in \curlyV_M} \v^\top \bSigmaw \v =  r \tilde{r} M^2 \frac{B n_B^2}{4(n_B-1)} =  r \tilde{r} M^2 \frac{n^2}{2n-B},
 \eneqn
 
\noindent where the last equality follows from the relation $n_B= 2n / B$.  Since $\frac{n^2}{2n-B}$ is monotone in $B$, it is minimized when $B$ is minimal, i.e., when $B=1$, which corresponds to i/BCRD.

In order to show Equation~\ref{eq:key}, let ${\B}$ be the matrix $\frac{n_B}{n_B - 1} \I_{n_B} - \frac{1}{n_B - 1} \J_{n_B}$ where $\J$ denotes the matrix of all ones. Recall that under $BL(B)$,  $\bSigmaw$ is block diagonal $B$ blocks of the form $(1-b){\B}$; see Appendix \ref{app:ibcrd_minimax}. Therefore, for block designs
\[
\v^\top \bSigmaw \v= r \tilde{r} \sum_{j =1}^B \v_j^\top \B \v_j, 
\]
where the $\v_j$'s are sub-vectors of $\v$. In order to maximize the quadratic form, one can maximize each term $\v_j^\top \B \v_j$ separately. By Lemma~\ref{lemm:maximal_vector}, the maximizer of each term contains  
$n_b/2$ entries with $M$, and $n_b/2$ entries with 0 (assuming $n_B$ is even). Equation~\ref{eq:etaB} within the lemma implies Equation~\ref{eq:key}, completing the proof of Theorem~\ref{thm:ibcrd_minimax_blocks}.

\begin{lemma}\label{lemm:maximal_vector}
The vector that maximizes the quadratic form ${\boldsymbol \eta}^\top {\bf B} {\boldsymbol \eta}$ for ${\boldsymbol \eta} \in \curlyV^{n_B}$ such that $\min(\eta_1,\ldots,\eta_{n_B}) \ge 0$  and $\max(\eta_1,\ldots,\eta_{n_B}) \le M$ when $n_B$ is even, has $n_B/2$ entries with $M$, and $n_B/2$ entries with 0.

\begin{proof}
First, since ${\boldsymbol \eta}^\top {\bf B} {\boldsymbol \eta}$ is a convex function in ${\boldsymbol \eta}$,  every entry of the maximizer must contain  either 0 or $M$. 
Also, by the definition of ${\B}$, ${\boldsymbol \eta}^\top {\bf B} {\boldsymbol \eta}$ is invariant under permutations of the entries of the vector ${\boldsymbol \eta}$. Thus, it is enough to consider vectors of the form  ${\boldsymbol \eta}_k$, which is defined to be a vector whose first $k$ entries are $M$ and the rest are zero. We have

\bneqn\label{eq:etaB}
{\boldsymbol \eta}_k^\top {\bf B} {\boldsymbol \eta}_k=
{\boldsymbol \eta}_k^\top \left( \frac{n_B}{n_B - 1} \I_{n_B} - \frac{1}{n_B - 1} \J_{n_B}\right) {\boldsymbol \eta}_k=\frac{M^2}{n_B -1}(n_Bk -k^2). 
\eneqn

\noindent Therefore, 

\beqn\label{eq:etaB1}
{\boldsymbol \eta}^\top_{k+1} {\bf B} {\boldsymbol \eta}_{k+1}-{\boldsymbol \eta}_{k}^\top {\bf B} {\boldsymbol \eta}_{k}= 
\frac{M^2}{n_B -1}[ n_B- (2k-1)].
\eeqn

\noindent It follows that if  $k < n_B/2$, then ${\boldsymbol \eta}^\top_{k+1} {\bf B} {\boldsymbol \eta}_{k+1}-{\boldsymbol \eta}_{k}^\top {\bf B} {\boldsymbol \eta}_{k}$ is non-negative, and if $k > n_B/2$, the difference is non-positive. Therefore, the maximizer is when $k=n_B/2$.
\end{proof}
\end{lemma}

\subsection{Expectation of the MSE}\label{app:expe_mse}

For this section and following sections, let $\Z := \Z_T / r + \Z_C / \rcontrol$ and $\muvec := \muvec_T  / r + \muvec_C  / \rcontrol$. 

\beqn
\expesub{\Z}{\msesub{\W}{\tauhat}} &=& \expesub{\Z_T, \Z_C}{\oneover{4} 
\parens{
    \frac{\muvec_T + \Z_T}{n_T} + \frac{\muvec_C + \Z_C}{n_C}
    }^\top \bSigmaw \parens{
    \frac{\muvec_T + \Z_T}{n_T} + \frac{\muvec_C + \Z_C}{n_C}
    }
}\\
&=& \oneover{4n^2} \expesub{\Z}{ 
\parens{
    \muvec + \Z
    }^\top \bSigmaw \parens{
    \muvec + \Z
    }
}\\
&=& \oneover{4n^2} \expesub{\Z}{ 
  \muvec^\top \bSigmaw \muvec +
  2\muvec^\top \bSigmaw \Z +
  \Z^\top \bSigmaw \Z
}\\
&=& \oneover{4n^2} \parens{\muvec^\top \bSigmaw \muvec + \expesub{\Z}{
  \Z^\top \bSigmaw \Z
}}\\
&=& \oneover{4n^2} \parens{\muvec^\top \bSigmaw \muvec + \tr{\bSigmaw\bSigmaZ}}\\
\eeqn

\noindent where $\bSigmaZ$ is the diagonal variance-covariance matrix of $\Z$. To compute the trace, we need to consider only the diagonal entries of $\bSigmaw \bSigmaZ$. The $i$th diagonal entry of this product is ${\bv{\Sigma}_{W}}_{i\, \cdot} {\bv{\Sigma}_{Z}}_{\cdot\, i} = {\bv{\Sigma}_{W}}_{i\, i} {\bv{\Sigma}_{Z}}_{i\, i}$ due to $\bSigmaZ$ being diagonal. From Equation~\ref{eq:diagonal_of_bSigmaW}, we know that all the diagonal entries of $\bSigmaw$ are $r\rcontrol$ and thus 

\bneqn\label{eq:trace_equality}
\tr{\bSigmaw \bSigmaZ} = r\rcontrol\, \tr{\bSigmaZ}.
\eneqn

\noindent Upon substitution, the result is

\beqn\footnotesize
\expesub{\Z}{\msesub{\W}{\tauhat}} &=& \oneover{4}\parens{
    \frac{\muvec_T}{n_T} + \frac{\muvec_C}{n_C}}^\top
\bSigmaw \parens{\frac{\muvec_T}{n_T} + \frac{\muvec_C}{n_C}
} 
+ 
\underbrace{\frac{1}{4n^2} r\rcontrol\,  \tr{\bSigmaZ}} \\
&=& \oneover{4n^2}\parens{
    \frac{\muvec_T}{r} + \frac{\muvec_C}{\rcontrol}}^\top
\bSigmaw \parens{\frac{\muvec_T}{r} + \frac{\muvec_C}{\rcontrol}
} 
+ c_Z.
\eeqn

\noindent where the underbraced term is design-independent and denoted as $c_Z$.

\subsection{Minimal block-size design is optimal for the mean criterion}\label{app:opt_expe_mse}

 Since $k$ is the minimal possible block-size, it follows that $m$ is divisible by $k$. Also, since $\bSigmaw$ is block-diagonal under block-designs, it it enough to consider sub-vectors of size $m$. Let $\muvec_m=(\mu_1,\ldots,\mu_m)^\top$ where $\mu_1 \le \mu_2 \le \cdots \le \mu_m$ and define $m_k=m/k$. Let $\vec{\muvec}_j$ denote  the sub-vector of the $j$-th block, and let $SSE(\vec{\muvec}_j)= \| \muvec_j -\onevec \bar{\muvec}_j \|^2$. 
We need to show that
\[
\muvec_m^\top \bSigma_{BL(m_k)} \muvec_m = (1-b) \sum_{j=1}^{m_k} \vec{\muvec}_j^\top \left(\frac{k}{k - 1} \I_{k} - \oneover{k - 1} \J_{k}\right)\vec{\muvec}_j=(1-b) \frac{k}{k - 1} SSE(\vec{\muvec}_j),
\]
is smaller than 
\[
\muvec_m^\top \bSigma_{BL(1)} \muvec_m = (1-b) \muvec_m^\top \left(\frac{m}{m - 1} \I_{m} - \oneover{m - 1} \J_{m}\right)\muvec_m=(1-b)\frac{m}{m-1} SSE(\muvec_m).
\]
That is, we need to show that 
\begin{equation} \label{eq:blocks}
  \frac{m}{m-1} SSE(\muvec_m)- \frac{k}{k - 1}\sum_{j=1}^{m_k}  SSE(\vec{\muvec}_j)  \ge 0.  
\end{equation}
We prove Inequality~\ref{eq:blocks} by induction on $m_k$ using Lemma \ref{lem:blocks} below. If $m_k=2$,  Inequality~\ref{eq:blocks} follows  
immediately from Lemma~\ref{lem:blocks}. 

Assume that  Inequality~\ref{eq:blocks} is true for $m_k=B$.  Consider $m_k=B+1$. Assume $\muvec_m$  is an ordered vector of length $m=k(B+1)$. By the induction assumption,
\beqn
\frac{k}{k-1}\sum_{j=1}^B SSE(\vec{\muvec}_j) \le \frac{kB}{kB-1}SSE(\vec{\muvec}_1 \cup \ldots \cup \vec{\muvec}_B).
\eeqn
We now apply Lemma \ref{lem:blocks} to the $m=k(B+1)$ elements by considering $\ell=Bk$. This implies that
\beqn
 \frac{kB}{kB-1}SSE(\vec{\muvec}_1 \cup \ldots \cup \vec{\muvec}_B)+\frac{k}{k-1}SSE(\vec{\muvec}_{B+1}) \le \frac{m}{m-1}SSE(\muvec_m).
\eeqn
Inequality~\ref{eq:blocks} follows immediately  by combining these two inequalities.

\begin{lemma}\label{lem:blocks}

For a vector ${\boldsymbol \nu}$ of length $l$ define $f({\boldsymbol \nu})=\frac{l}{l-1}SSE({\boldsymbol \nu})$. Consider a division of $\muvec_m$
into two sub-vectors $\underline{\muvec}_1,\underline{\muvec}_2$ of size $\ell$ and $m-\ell$, respectively, i.e., $\underline{\muvec}_1=(\mu_1,\ldots,\mu_\ell)$ and $\underline{\muvec}_2=(\mu_{\ell+1},\ldots,\mu_m)$. Then,
\begin{equation} \label{eq:lem}
   f(\muvec_m) - f(\underline{\muvec}_1) - f(\underline{\muvec}_2) \ge 0.
\end{equation}

\begin{proof}
Since we can  subtract from the entries $\mu_i$ a constant without changing the values of $f(\underline{\muvec}_1)$, $f(\underline{\muvec}_2)$ and $f(\muvec)$, we can assume without loss of generality that the entries in $\underline{\muvec}_1$ are non-positive and the entries in $\underline{\muvec}_2$, non-negative.

Define a vector $\tilde{\muvec}$ of length $m$ as follows 
\[
\tilde{\mu}_i= \left\{ \begin{array}{cc} 0 & i=\ell+1\\ \mu_m+\mu_{\ell+1} & i=m\\ \mu_i & \text{otherwise} \end{array}  \right. .
\]
and let $\underline{\tilde{\muvec}}_1$ and $\underline{\tilde{\muvec}}_2$ be the corresponding sub-vectors.  Notice that $\underline{\tilde{\muvec}}_1=\underline{{\muvec}}_1$ and $\bar{\underline{\tilde{\muvec}}}_2=\bar{\underline{{\muvec}}}_2$ and $\bar{{\tilde{\muvec}}}=\bar{{{\muvec}}}$. A trivial calculation shows that 
\[
SSE( \underline{\tilde{\muvec}}_2 ) - SSE(\underline{{\muvec}}_2)=
SSE(\tilde{\muvec})-SSE(\muvec_m)= 2\mu_m \mu_{\ell+1} >0.
\]
Thus, when replacing $\muvec_m$ with $\tilde{\muvec}$,  $SSE({{\muvec}}_2)$ goes up by the same amount as $SSE(\muvec_m)$.
But, since $\frac{m-\ell}{m-\ell-1}>\frac{m}{m-1}$, 
\[
f(\muvec_m)- f(\underline{\muvec}_2) >  f(\tilde{\muvec}) - f(\underline{\tilde{\muvec}}_2).
\]
It follows that replacing $\muvec_m$ with $\tilde{\muvec}$ reduces the right hand-side of Inequality~\ref{eq:lem} and therefore we can assume that $\mu_{\ell+1}=0$. 
The same argument shows that we can also assume that $\mu_{\ell+2}=0$. Continuing in the same manner implies  that we can consider $\underline{\muvec}_2=(0,\ldots,0,c)$ for $c\ge 0$ as the worst case. A symmetric argument shows that $\underline{\muvec}_1=(a,0,\ldots,0)$ were $a \le 0$ is the worst case for $\underline{\muvec}_1$. 

It follows that it is enough to prove Inequality~\ref{eq:lem} for $\muvec_m=(a,0,\ldots,0,c)$  where $a\le0$ and $c\ge0$. First consider $\underline{\muvec}_2$. Since the average of the $m-\ell$ elements in $\underline{\muvec}_2$ is $\frac{c}{m-\ell}$, that implies
\beqn
SSE(\underline{\muvec}_2)&=&(m-\ell-1)\frac{c^2}{(m-\ell)^2}+\left(c-\frac{c}{m-\ell}\right)^2\\
&=&\frac{c^2}{(m-\ell)^2}(m-\ell-1+(m-\ell-1)^2)\cr
&=&c^2\frac{m-\ell-1}{m-\ell}.
\eeqn
Therefore, $f(\underline{\muvec}_2)=c^2$. A parallel argument shows that $f(\underline{\muvec}_1)=a^2$. 
We now need $SSE(\muvec_m)$.
We have that $\bar{\muvec}_m=\frac{a+c}{m}$ and hence 
\beqn
SSE(\muvec_m)&=&\left(a-\frac{a+c}{m}\right)^2+ (m-2)\left(\frac{a+c}{m}\right)^2+\left(c-\frac{a+c}{m}\right)^2\\
&=&a^2+c^2-m\left(\frac{a+c}{m}\right)^2=a^2+c^2-\frac{a^2+c^2}{m}-\frac{2ac}{m}\\
&=&\frac{m-1}{m}(a^2+c^2)-\frac{2ac}{m} \ge \frac{m-1}{m}(a^2+c^2) .
\eeqn
The last inequality follows because $c \ge 0$ and $a \le 0$. Therefore, $f(\muvec_m) \ge f(\underline{\muvec}_1)+ f(\underline{\muvec}_2)$ and Inequality~\ref{eq:lem} follows. 
\end{proof}
\end{lemma}

\subsection{Variance of the MSE}\label{app:var_mse}

We compute $\varsub{\Z}{\msesub{\W}{\tauhat}}$ below.

\bneqn\label{eq:var_mse}
&=& \varsub{\Z}{\oneover{4n^2} (\muvec + \Z)^\top \bSigmaw (\muvec + \Z)} \nonumber \\
&=& \oneover{16n^4}\varsub{\Z}{\muvec^\top \bSigmaw \muvec + 2\muvec^\top \bSigmaw \Z + \Z^\top \bSigmaw \Z} \nonumber \\
&=& \oneover{16n^4}\varsub{\Z}{2\muvec^\top \bSigmaw \Z + \Z^\top \bSigmaw \Z} \nonumber \\
&=& \oneover{16n^4} \parens{4\,\varsub{\Z}{\Z^\top \bSigmaw \muvec} + \varsub{\Z}{\Z^\top \bSigmaw \Z} + 4 \,\cov{\muvec^\top \bSigmaw \Z}{\Z^\top \bSigmaw \Z}} \nonumber \\
&=& \oneover{16n^4} \parens{ 4(\bSigmaw\muvec)^\top \bSigmaZ (\bSigmaw\muvec) + \varsub{\Z}{\Z^\top \bSigmaw \Z} + 4\,\cov{\muvec^\top \bSigmaw \Z}{\Z^\top \bSigmaw \Z} } \nonumber \\
&=& \oneover{16n^4} \parens{ 4\parens{\muvec^\top \bSigmaw \bSigmaZ \bSigmaw\muvec + \cov{\muvec^\top \bSigmaw \Z}{\Z^\top \bSigmaw \Z}} + \varsub{\Z}{\Z^\top \bSigmaw \Z} }
\eneqn

\noindent The covariance term is the expectation of the product of the two inputs minus the product of the expectations of the two inputs. Since $\expe{\Z} = \zerovec_{2n}$, the product of the expectations becomes zero. Hence,

\beqn
\cov{\muvec^\top \bSigmaw \Z}{\Z^\top \bSigmaw \Z} = \muvec^\top \bSigmaw \expesub{\Z}{\Z \Z^\top \bSigmaw \Z}
\eeqn

\noindent The $k$th element of the expectation of the cubic form above looks like

\beqn
\expesub{\Z}{Z_k \Z^\top \bSigmaw \Z} = \expesub{\Z}{Z_k \sum_{i,j} Z_i Z_j \sigmawij{i}{k}}
\eeqn

\noindent where $\sigmawij{i}{j}$ denotes the $i,j$th entry of the $\bSigmaw$ matrix. Due to the independence of the $Z_i$'s the terms in the above sum are zero except when $i=j=k$ which results in

\beqn
\expesub{\Z}{Z_k \Z^\top \bSigmaw \Z} = \expesub{\Z}{\sigmawij{k}{k} Z_k^3}.
\eeqn

\noindent Since the diagonal entries of $\bSigmaw$ are identically $r\rcontrol$ (see Equation~\ref{eq:diagonal_of_bSigmaW}), then the vector whose $k$th entry we were looking at is merely

\beqn
\expesub{\Z}{\Z \Z^\top \bSigmaw \Z} = r\rcontrol\gammaZ
\eeqn

\noindent where $\gammaZ$ is the vector whose elements are the third moments of $\Z$. Thus,

\bneqn\label{eq:cov_term}
\cov{\muvec^\top \bSigmaw \Z}{\Z^\top \bSigmaw \Z} = r\rcontrol\muvec^\top \bSigmaw \gammaZ.
\eneqn

\noindent We now focus on the remaining variance term in the overall variance which can be expressed as

\bneqn\label{eq:variance_to_quartic}
\varsub{\Z}{\Z^\top \bSigmaw \Z} &=& \expesub{\Z}{(\Z^\top \bSigmaw \Z)^2} -  (\expesub{\Z}{\Z^\top \bSigmaw \Z})^2 \nonumber \\
&=& \expesub{\Z}{(\Z^\top \bSigmaw \Z)^2} - r^2 \rcontrol^2 \tr{\bSigmaZ}^2
\eneqn

\noindent where the last equality follows from Equation~\ref{eq:trace_equality}. We now compute  $\expesub{\Z}{(\Z^\top \bSigmaw \Z)^2}$, the quartic form in the above expression. We define $\sigmazi{i} := \expe{Z_i^2}$ which is also the $i,i$th entry of the $\bSigmaZ$ matrix. The quartic form can be expanded in scalar form as follows:

\bneqn\label{eq:quartic_expansion}
\expe{(\Z^\top \bSigmaw \Z)^2} = \sum_{i,j,k,\ell \in \braces{1, \ldots, 2n}} \sigmawij{i}{j} \sigmawij{k}{\ell} \expe{Z_i Z_j Z_k Z_\ell}
\eneqn

\noindent where the expectation operator is taken over the relevant $Z_i$'s. Since $\expe{Z_i} = 0$ for all subjects, the expectation above simplifies to

\beqn
\expe{Z_i Z_j Z_k Z_\ell} = \begin{cases}
\expe{Z_i^4} \quad \text{if} \quad i=j=k=\ell \\
\sigmazi{i}\sigmazi{k} ~~~\quad \text{if} \quad i=j \neq k=\ell \\
\sigmazi{i}\sigmazi{j} ~~~\quad \text{if} \quad i=k \neq j=\ell ~~ \text{or} ~~ i=\ell \neq j=k \\
0 ~~~~~~~\quad \text{otherwise}
\end{cases}
\eeqn

\noindent Using the above result, we write the quartic expansion of Equation~\ref{eq:quartic_expansion} as

\beqn
= \sum_{i=1}^{2n} \sigmawij{i}{i}^2 \expe{Z_i^4}  + 
\sum_{i=j \neq k=\ell} \sigmawij{i}{i}\sigmawij{k}{k} \sigmazi{i}\sigmazi{k} + 
\sum_{i=k \neq j=\ell} \sigmawij{i}{j} \sigmawij{i}{j} \sigmazi{i}\sigmazi{j} + 
\sum_{i=\ell \neq j=k} \sigmawij{i}{j} \sigmawij{j}{i} \sigmazi{i}\sigmazi{j}.
\eeqn

\noindent Given that $\bSigmaw$ is symmetric with a diagonal of entries that are identically $r\rcontrol$ (see Equation~\ref{eq:diagonal_of_bSigmaW}), the above simplifies to

\bneqn\label{eq:quartic_expansion_simplified}
&=& \r^2 \rcontrol^2 \sum_{i=1}^{2n} \expe{Z_i^4}  + 
r^2 \rcontrol^2 \sum_{i=j \neq k=\ell} \sigmazi{i}\sigmazi{k} + 
\sum_{i=k \neq j=\ell} \sigmawij{i}{j}^2 \sigmazi{i}\sigmazi{j} + 
\sum_{i=\ell \neq j=k} \sigmawij{i}{j}^2 \sigmazi{i}\sigmazi{j} \nonumber \\
&=& r^2 \rcontrol^2 \parens{
    \sum_{i=1}^{2n} \expe{Z_i^4}  + 
    \sum_{i \neq j} \sigmazi{i}\sigmazi{j} 
} + 
2\sum_{i \neq j} \sigmawij{i}{j}^2 \sigmazi{i}\sigmazi{j}.
\eneqn

\noindent Consider the following matrix and its entries:

\beqn
(\bSigmaw\bSigmaZ)_{i,j} = \sigmawij{i}{j} \sigmazi{j}
\eeqn

\noindent We now square this matrix and examine only the result's diagonal entries

\beqn
[(\bSigmaw\bSigmaZ)^2]_{i,i} = \bracks{\sigmawij{i}{1} \sigmazi{1} \ldots \sigmawij{i}{2n} \sigmazi{2n}} \bracks{\sigmawij{1}{j} \sigmazi{j} \ldots \sigmawij{2n}{j} \sigmazi{j}}^\top = \sum_{j=1}^{2n} (\sigmawij{i}{j} \sigmazi{j}) (\sigmawij{j}{i} \sigmazi{i}) = \sum_{j=1}^{2n} \sigmawij{i}{j}^2 \sigmazi{j} \sigmazi{i}
\eeqn

\noindent Thus,

\bneqn\label{eq:trSigma} \footnotesize
2\,\tr{(\bSigmaw\bSigmaZ)^2} = 2 \sum_{i=1}^{2n} \sum_{j=1}^{2n} \sigmawij{i}{j}^2 \sigmazi{j} \sigmazi{i} =2 \sum_{i\neq j} \sigmawij{i}{j}^2 \sigmazi{j} \sigmazi{i} + 2 \sum_{i = j} \sigmawij{i}{j}^2 \sigmazi{j} \sigmazi{i} = 2 \sum_{i\neq j} \sigmawij{i}{j}^2 \sigmazi{j} \sigmazi{i} + 2 r^2 \rcontrol^2 \sum_{i=1}^{2n} \sigmazi{i}^2
\eneqn

\noindent Since the above shares the first term with the simplified quartic form of Equation~\ref{eq:quartic_expansion_simplified}, we now write the quartic form as a function of this trace,

\beqn
\expe{(\Z^\top \bSigmaw \Z)^2} = r^2 \rcontrol^2 \parens{
    \sum_{i=1}^{2n} \expe{Z_i^4}  + 
    \sum_{i \neq j} \sigmazi{i}\sigmazi{j} - 
    2 \sum_{i=1}^{2n} \sigmazi{i}^2
} + 
2\,\tr{(\bSigmaw\bSigmaZ)^2}.
\eeqn

\noindent We now substitute this into the variance expression of Equation~\ref{eq:variance_to_quartic} to obtain

\beqn
\varsub{\Z}{\Z^\top \bSigmaw \Z} = r^2 \rcontrol^2 \parens{
    \sum_{i=1}^{2n} \expe{Z_i^4}  + 
    \sum_{i \neq j} \sigmazi{i}\sigmazi{j} - 
    2 \sum_{i=1}^{2n} \sigmazi{i}^2 -
    \tr{\bSigmaZ}^2
} + 
2\,\tr{(\bSigmaw\bSigmaZ)^2}.
\eeqn

\noindent Substituting the expansion of $\tr{\bSigmaZ}^2 = \sum_{i \neq j} \sigmazi{i}\sigmazi{j} + \sum_{i=1}^{2n} \sigmazi{i}^2$ we arrive at

\bneqn\label{eq:final_var_z_term}
\varsub{\Z}{\Z^\top \bSigmaw \Z} &=& r^2 \rcontrol^2 \parens{
    \sum_{i=1}^{2n} \expe{Z_i^4} - 
    3 \sum_{i=1}^{2n} \sigmazi{i}^2
} + 
2\,\tr{(\bSigmaw\bSigmaZ)^2} \nonumber \\
&=& r^2 \rcontrol^2 \kappa_Z + 2\,\tr{(\bSigmaw\bSigmaZ)^2}.
\eneqn

\noindent The terms $\expe{Z_i^4}- 3 \sigmazi{i}^2$ are a constant with respect to the design and we denote the sum of these terms by $\kappa_Z$. We chose this notation as the quantity expressed by the sum of these terms is similar to the sum of the excess kurtosis of the $Z_i$'s. 

Substituting the two intermediate results of Equations~\ref{eq:cov_term} and \ref{eq:final_var_z_term} into the overall variance expression of Equation~\ref{eq:var_mse}, we arrive at the final expression

\beqn
\varsub{\Z}{\msesub{\W}{\tauhat}} = \oneover{16n^4}\parens{ 4\parens{
  \muvec^\top \bSigmaw \bSigmaZ \bSigmaw\muvec + r \rcontrol \muvec^\top \bSigmaw \gammaZ
} + r^2 \rcontrol^2 \kappa_Z + 2\,\tr{(\bSigmaw\bSigmaZ)^2}
}.
\eeqn

\subsection{An asymptotic lower bound for the tail criterion (Proof of Theorem \ref{thm:lower})}\label{app:lower}

If $\mathcal{B}_1$ does not converge, then the following argument holds for any convergent subsequence of $\mathcal{B}_1$. 

Recall the definition of $\tilde{Q}$ in Equation~\ref{eq:Qtilde}. If $\mathcal{B}_1$ is of order $n$, then $\lim_{n\to \infty}\frac{1}{\sqrt n} \tilde{Q} = \infty$ and  the lower bound (Equation~\ref{eq:lower}) clearly holds. So, we now assume that $\mathcal{B}_1=o(n)$. In this case, by Equation~\ref{eq:Cbound} below, $|\mathcal{S}|=o(n)$ and since $\mathcal{R}$ is at-least of order $n$ (see below), then $\mathcal{S}$ is negligible. The terms $\mathcal{B}_1$ and $\mathcal{B}_2$ are positive. Therefore,

\beqn
\liminf_{n \to \infty} \frac{1}{\sqrt n} \tilde{Q} &=& \liminf_{n \to \infty} \frac{1}{\sqrt n} \left(\mathcal{B}_1 + c \times \sqrt{r^2 \tilde{r}^2 \kappa_Z + 4\parens{\mathcal{B}_2 + \mathcal{S}} + 2 \mathcal{R} } \right)\\
& \ge & \liminf_{n \to \infty}   c \times \sqrt{ r^2 \tilde{r}^2\kappa_Z/n  + 2 \mathcal{R}/n }.
\eeqn

\noindent Now, by Assumption~\ref{ass:boundedness3}, $\kappa_Z/n \to 2\bar{\kappa}_Z$, and by Equation~\ref{eq:trSigma}

\[
\frac{ \mathcal{R}}{n} = \frac{  \tr{(\bSigmaw\bSigmaZ)^2}}{n} \ge \frac{ r^2 \tilde{r}^2 \sum_{i=1}^{2n} \rho_i^2}{n} = r^2 \tilde{r}^2 2 \frac{1}{2n} \| \brho\|^2 \to  2r^2 \tilde{r}^2 \bar{\brho},
\]

\noindent where the last equality is due to Assumption \ref{ass:boundedness1}. This completed the proof of Theorem \ref{thm:lower}.

\subsection{Optimal Design is Harmony (Proof of Theorem \ref{thm:H})}\label{app:harmony}

The arguments of Sections~\ref{subsubsec:tail_PB} and \ref{subsubsec:tail_BCRD} show that under both PB and i/BCRD, $\tilde{Q}$ (Equation~\ref{eq:Qtilde}) is of order $n$ and therefore, $\frac{1}{\sqrt{n}} \tilde{Q} \to \infty$.


The only thing that is left to show is that under PB and under equal allocation, $\mathcal{R}$ is of order $n^2$. Since PB is composed from a single vector $\w_*$ and its mirror $-\w_*$, we have that under PB, $\bSigmaw=\w_* \w_*^\top$.
Recall that $\brho$ denotes the diagonal of $\bSigmaZ$ and let $\brho_{\w_*}$ denote the element-wise product of $\brho$ and ${\w_*}$.  
Therefore, under PB,

\beqn
\mathcal{R} &=& \tr{(\bSigmaw\bSigmaZ)^2}= \tr{ \w_* \w_*^\top \bSigmaZ \w_* \w_*^\top \bSigmaZ  }=
\tr{ \w_* \brho_{\w_*}^\top \w_* \brho_{\w_*}^\top } \\
&=& (\brho_{\w_*}^\top \w_*)^2 = \left(\sum_{i=1}^n \w_{*,i}^2 \rho_i \right)^2
= \left(\sum_{i=1}^n  \rho_i \right)^2, 
\eeqn
which is of order $n^2$ by Assumption \ref{ass:boundedness1}.

\subsection{Asymptotics of the tail criterion (Proof of Theorem \ref{thm:Q})}\label{app:asymp}




We proceed by bounding $\mathcal{B}_2$ and $\mathcal{S}$ and then computing term $\mathcal{R}$. By Assumption \ref{ass:boundedness1},

\bneqn\label{eq:Bbound}
\mathcal{B}_2= \muvec^\top \bSigmaw \bSigmaZ \bSigmaw \muvec \le \bar{C}_{\brho} \| \bSigmaw \muvec \|^2
= \bar{C}_{\brho} r \tilde{r} \sum_{b =1}^B \| \B \muvec_b \|^2
= \bar{C}_{\brho} r \tilde{r} \frac{n_B}{n_B-1} \mathcal{B}_1 \le 2 \bar{C}_{\brho} r \tilde{r} \mathcal{B}_1.
\eneqn

\noindent By the Cauchy-Schwarz inequality,

\bneqn\label{eq:Cbound}
|\mathcal{S}| &=& |\muvec^\top \bSigmaw \gammaZ| \le \sqrt{\muvec^\top \bSigmaw \muvec} \sqrt{\gammaZ^\top \bSigmaw \gammaZ} = \sqrt{\mathcal{B}_1 } \sqrt{\gammaZ^\top \bSigmaw \gammaZ} \nonumber \\
&\le& \sqrt{\mathcal{B}_1  \lambda_{max} (\bSigmaw)} \| \gammaZ \| \le   2 r \tilde{r} \sqrt{  C_{\boldsymbol{\gamma}} n \mathcal{B}_1  },  
\eneqn 

\noindent where the last inequality is true because for block designs $\lambda_{max} (\bSigmaw)= r \tilde{r} \frac{n_B}{n_B-1} \le 2 r \tilde{r} $ and due to Assumption \ref{ass:boundedness2}.

Let $\brho_b^2$ be the vector whose entries is the entries of $\brho_b$ squared. We have

\beqn
\mathcal{R} &=& \tr{(\bSigmaw\bSigmaZ)^2}=r^2 \tilde{r}^2\sum_{b=1}^B \tr{ (\B {\rm diag}(\brho_b ))^2 }
\\
&=  &\frac{r^2 \tilde{r}^2}{ (n_B - 1)^2} \sum_{b=1}^B
\tr{ (n_B {\rm diag}(\brho_b ) - \J_{n_B} {\rm diag}(\brho_b ) )^2 }\\
& = & \frac{r^2 \tilde{r}^2}{(n_B - 1)^2} \sum_{b=1}^B \left\{ n_B^2 \tr{  {\rm diag}(\brho_b^2)} - 2n_B \tr {\J_{n_B} {\rm diag}(\brho_b^2 )} + \tr{ [\J_{n_B} {\rm diag}(\brho_b )]^2}  \right\}.
\eeqn

\noindent We now compute the above traces:

\beqn
\tr{  {\rm diag}(\brho_b^2 )}&= & \| \brho_b \|^2;
\\
\tr {\J_{n_B} {\rm diag}(\brho_b^2 )} &=  & \tr { \onevec_{n_B} \onevec_{n_B} ^T {\rm diag}(\brho_b^2 ) } = \tr { \onevec_{n_B} (\brho_b^2 )^T }= \| \brho_b \|^2 ;  
\\
\tr{ [\J_{n_B} {\rm diag}(\brho_b )]^2}&=& \tr{  \onevec_{n_B} \onevec_{n_B} ^T {\rm diag}(\brho_b ) \onevec_{n_B} \onevec_{n_B} ^T {\rm diag}(\brho_b )} \\
&=& \tr{  \onevec_{n_B} \brho_b ^T   \onevec_{n_B} \brho_b ^T }= \brho_b ^T   \onevec_{n_B} \tr{  \onevec_{n_B}  \brho_b ^T }=  (\brho_b ^T   \onevec_{n_B})^2.
\eeqn

\noindent Putting all terms together we have

\beqn
\mathcal{R}&=& \frac{r^2 \tilde{r}^2}{(n_B - 1)^2} \sum_{b=1}^B \left\{ n^2_B  \| \brho_b \|^2 - 2n_B \| \brho_b \|^2 + (\brho_b ^T   \onevec_{n_B})^2   \right\}\\
&=& \frac{r^2 \tilde{r}^2}{(n_B - 1)^2} \sum_{b=1}^B \left\{ n_B \| \brho_b \|^2 (n_B-1) + (\brho_b ^T   \onevec_{n_B})^2 - n_B \| \brho_b \|^2  \right\} \\
&=& \frac{ r^2 \tilde{r}^2 n_B \|\brho\|^2  }{n_B - 1}+ \frac{r^2 \tilde{r}^2}{(n_B - 1)^2} \sum_{b=1}^B \left\{  (\brho_b ^T   \onevec_{n_B})^2 - n_B \| \brho_b \|^2  \right\}.
\eeqn

\noindent Now notice that

\beqn
\brho_b^T \B \brho_b =  \brho_b^T \left( \frac{n_B}{n_B - 1} \I_{n_B} - \frac{1}{n_B - 1} \J_{n_B} \right) \brho_b = \frac{1}{n_B -1} \left( n_B \| \brho_b \|^2 - (\brho_b ^T   \onevec_{n_B})^2 \right) .
\eeqn

\noindent Therefore,

\bneqn\label{eq:Dcomputation}
\mathcal{R} &=&  r^2 \tilde{r}^2\frac{n_B \|\brho\|^2 - \sum_{b=1}^B \brho_b^T \B \brho_b  }{n_B -1 } \nonumber \\
&=&  r^2 \tilde{r}^2 \frac{4n^2 \frac{1}{2n} \|\brho\|^2 - B \sum_{b=1}^B \brho_b^T \B \brho_b  }{2n -B } \nonumber \\
&=& \frac{  r^2 \tilde{r}^2 4n^2 \frac{1}{2n} \|\brho\|^2 - r \tilde{r} B \brho^T \bSigmaw \brho }{2n -B }.
\eneqn

\noindent We now consider four regimes for $B$ which become the four results of Theorem~\ref{thm:Q}.

\subsection*{(i) $\frac{B}{n^{1/4}} \to 0$}

By Assumption \ref{ass:blocking1},

\beqn
\frac{1}{\sqrt{n}}{\tilde{Q}} > \frac{1}{\sqrt{n}} \mathcal{B}_1 \approx \frac{1}{\sqrt{n}} n/B^2 \to \infty. 
\eeqn

\subsection*{(ii) $\frac{B}{n^{1/4}} \to 1 $} 

By Assumption \ref{ass:blocking1} that in this case  $\frac{1}{\sqrt{n}} \mathcal{B}_1 \to C_{\mathcal{B}_1}$. By Assumption \ref{ass:boundedness3}, $\frac{\kappa_Z}{n} \to 2 \bar{\kappa}_Z$. The bounds for terms $\mathcal{B}_2$ and $\mathcal{S}$ (Equations~\ref{eq:Bbound} and \ref{eq:Cbound}) imply that $\sqrt{\frac{\mathcal{B}_2}{n}} \to 0$ and $\sqrt{\frac{\mathcal{S}}{n}} \to 0$ and hence those terms are negligible. For term $\mathcal{R}$, by Equation~\ref{eq:Dcomputation} and Assumptions \ref{ass:boundedness1} and \ref{ass:blocking2},

\beqn
\frac{\mathcal{R}}{n} =  r^2 \tilde{r}^2 \frac{4n \frac{1}{2n} \|\brho\|^2 }{2n -B }+ o(1) \to  2 r^2 \tilde{r}^2 \bar{\brho}. 
\eeqn

\noindent In conclusion, we have that under this regime

\beqn
\frac{2}{\sqrt{n}}{\tilde{Q}} \to C_{\mathcal{B}_1} + c  r\tilde{r}\sqrt{ 2(  \bar{\kappa}_Z +\bar{\brho}) }. 
\eeqn

\subsection*{(iii) $\frac{B}{n^{1/4}} \to \infty$ and $B/n \to \infty$} 

We have by Assumption \ref{ass:blocking1} that in this case  $\frac{1}{\sqrt{n}} \mathcal{B}_1 \to 0$ and similarly to before terms $\mathcal{B}_2$ and $\mathcal{S}$ are negligible. Also, as in the previous case $\frac{\mathcal{R}}{n}   \to 2  r^2 \tilde{r}^2\bar{\brho}$. 

\subsection*{(iv) $B/n \to \alpha$}

As in the previous case, terms $\mathcal{B}_1$, $\mathcal{B}_2$ and $\mathcal{S}$ are negligible and for term $\mathcal{R}$ we have that

\beqn
\frac{\mathcal{R}}{n} =  r^2 \tilde{r}^2 \frac{2 \frac{1}{2n} \|\brho\|^2 }{1 -B/(2n) }+ o(1) \to 2  r^2 \tilde{r}^2 \bar{\brho}\frac{1}{1-\alpha/2}. 
\eeqn

\section{Additional Simulations}\label{app:additional_simualtions}

\begin{figure}[htp]
\centering
     \begin{subfigure}[b]{\textwidth}
         \centering
         \includegraphics[width=6in]{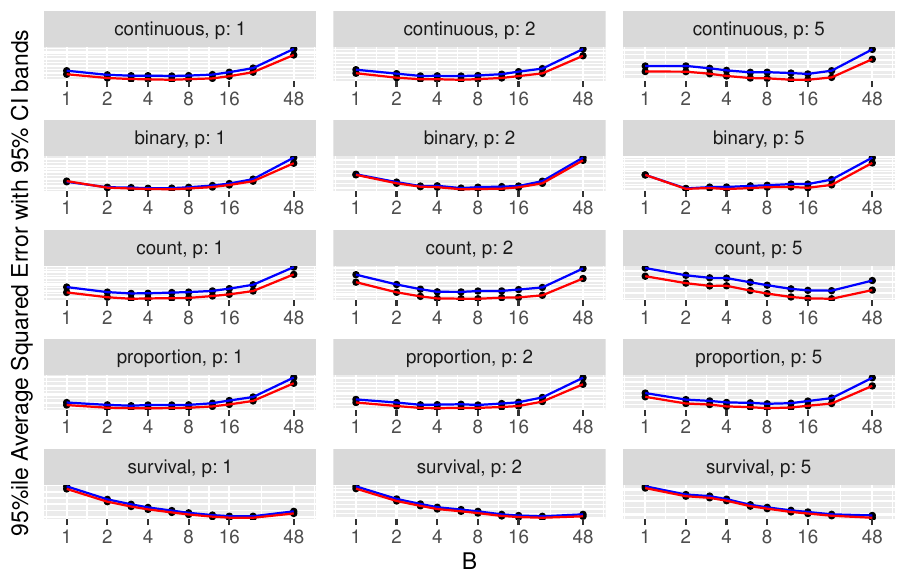}
         \caption{Equal Allocation}
         \label{fig:all_response_q_95_equal_allocation_b}
     \end{subfigure}\\
     \begin{subfigure}[b]{\textwidth}
         \centering
         \includegraphics[width=6in]{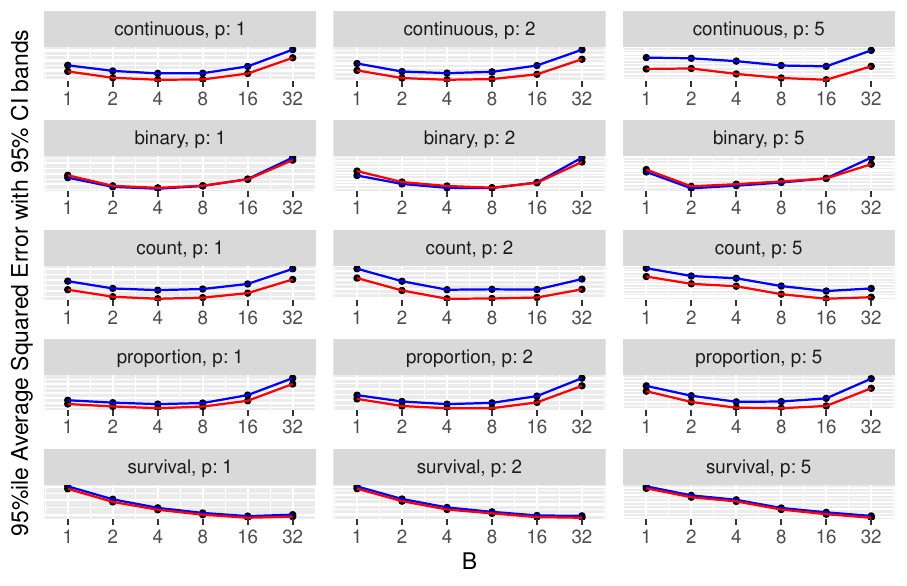}
         \caption{Unequal Allocation ($n_C:n_T$ = 2:1)}
         \label{fig:all_response_q_95_unequal_allocation_b}
     \end{subfigure} 
\caption{Simulation results for mean-centered exponential covariate realizations, $n=96$, $p = 1,2,5$, all response types and all appropriate number of blocks. The y-axis is relative to simulation and thus its values are unshown. Blue lines are empirical 95\%ile MSE Tails of (Equation~\ref{eq:actual_quantile_criterion}) and red lines are approximate quantiles computed via the average MSE plus 1.645 times the empirical standard error of the MSE (Equation~\ref{eq:approx_quantile_criterion}).}
\end{figure}

We run an identical simulation as in Section~\ref{sec:simulations} except we change the covariate distribution to be long-tailed (the mean-centered exponential) with parameter $1/\sqrt{3}$ for all response types except incidence which has parameter $\sqrt{3}$. These distributions were calibrated to have the same mean (zero) and variances as the covariate distribution found in the main text (cf. Table~\ref{tab:covariate_and_response_settings}). The results are found in Figures~\ref{fig:all_response_q_95_equal_allocation_b}-\ref{fig:all_response_q_95_unequal_allocation_b} and are qualitatively similar as the simulations found in the main text. This implies that our results and asymptotics are robust to the setting of long-tailed observed subject measurements.

       
       
       
       
       

\end{document}